\documentclass[11pt,english,reqno]{smfart}
\usepackage{times}\usepackage{mathptmx}
\def\baselinestretch{1.1}
\usepackage{amsmath}
\usepackage{amsthm}
\usepackage{amsfonts}
\usepackage{amssymb}

\newcommand{\BBC}[2]{
\Bigl( \!\!\! \begin{array}{c} #1 \\[-0.6mm] #2 
             \end{array} \!\!\! \Bigr)_b}

\theoremstyle{plain}
\newtheorem{thm}{Theorem}

\newtheorem{propn}{Proposition}
\newtheorem{lem}{Lemma}
\newtheorem{cor}{Corollary}

\newtheorem{defn}{Definition}

\theoremstyle{remark}
\newtheorem{rem}{Remark}

\renewcommand{\theequation}{\thesection.\arabic{equation}} 
\newcommand{\id}{{\rm id}}
\newcommand{\ot}{\otimes}
\newcommand{\ra}{\rightarrow}
\newcommand{\ti}{\times}
\newcommand{\fr}[2]{{\textstyle \frac{#1}{#2} }}

\newcommand{\fsl}{{\mathfrak s}{\mathfrak l}}

\newcommand{\bra}{\langle}
\newcommand{\ket}{\rangle}
\newcommand{\al}{\alpha}
\newcommand{\be}{\beta}
\newcommand{\ga}{\gamma}

\newcommand{\de}{\delta}
\newcommand{\De}{\Delta}
\newcommand{\ep}{\epsilon}

\newcommand{\om}{\omega}
\newcommand{\Om}{\Omega}
\newcommand{\si}{\sigma}

\newcommand{\CA}{{\mathcal A}}
\newcommand{\CB}{{\mathcal B}}

\newcommand{\CD}{{\mathcal D}}
\newcommand{\CH}{{\mathcal H}}

\newcommand{\CK}{{\mathcal K}}

\newcommand{\CO}{{\mathcal O}}  
\newcommand{\CP}{{\mathcal P}}  
\newcommand{\CS}{{\mathcal S}}
\newcommand{\CT}{{\mathcal T}}
\newcommand{\CU}{{\mathcal U}}

\newcommand{\SB}{{\mathsf B}}
\newcommand{\SC}{{\mathsf C}}

\newcommand{\SE}{{\mathsf E}}
\newcommand{\SF}{{\mathsf F}}

\newcommand{\SH}{{\mathsf H}}

\newcommand{\SJ}{{\mathsf J}}
\newcommand{\SK}{{\mathsf K}}

\newcommand{\SO}{{\mathsf O}}  
\newcommand{\SP}{{\mathsf P}}  
  
\newcommand{\SR}{{\mathsf R}}

\newcommand{\SX}{{\mathsf X}}
\newcommand{\SY}{{\mathsf Y}}

\newcommand{\se}{{\mathsf e}}
\renewcommand{\sf}{{\mathsf f}}
\newcommand{\sh}{{\mathsf h}}

\newcommand{\spp}{{\mathsf p}}

\renewcommand{\ss}{{\mathsf s}}

\newcommand{\sx}{{\mathsf x}}

\newcommand{\one}{{\mathfrak 1}}
\newcommand{\two}{{\mathfrak 2}}
\newcommand{\three}{{\mathfrak 3}}

\newcommand{\BR}{{\mathbb R}}

\newcommand{\BC}{{\mathbb C}}

\newcommand{\BZ}{{\mathbb Z}}

\newcommand{\CGC}[6]{\big[
\begin{smallmatrix} #1 & #3 & #5 \\
#2 & #4 & #6
\end{smallmatrix}\big]} 
\newcommand{\CGCb}[6]{\Big[
\begin{smallmatrix} #1 & #3 & #5 \\
#2 & #4 & #6
\end{smallmatrix}\Big]} 
\newcommand{\CGCtr}[5]{C_{#1}{\big[
\begin{smallmatrix} #2 & #4  \\
#3 & #5 
\end{smallmatrix}\big]}}

\DeclareMathOperator*{\Res}{Res}
\renewcommand{\=}[1]{\stackrel{(\ref{#1})}{=}} 
\newcommand{\rf}[1]{(\ref{#1})}

\begin{document}
\title[R-operator, co-product and Haar-measure for the modular double]
{R-operator, co-product and Haar-measure for the modular double
of $U_q(\fsl(2,\BR))$}
\author{A.G. Bytsko 
  \ and \ J. Teschner}
\address{Institut f\"ur theoretische Physik\\
Freie Universit\"at Berlin\\
Arnimallee 14, 14195 Berlin, Germany}
\email{bytsko@physik.fu-berlin.de, teschner@physik.fu-berlin.de}
\begin{abstract}
A certain class of unitary 
representations of $U_q(\fsl(2,\BR))$ has the property of being
simultanenously a representation of 
$U_{\tilde{q}}(\fsl(2,\BR))$ for a particular choice of $\tilde{q}(q)$. 
Faddeev has 
proposed to unify the quantum groups $U_q(\fsl(2,\BR))$ and 
$U_{\tilde{q}}(\fsl(2,\BR))$ into some enlarged object for which he has
coined the name ``modular double''. 

We study the R-operator, the co-product and the 
Haar-measure for the modular double
of $U_q(\fsl(2,\BR))$ and establish their main properties.
In particular it is shown that the Clebsch-Gordan maps constructed
in \cite{PT} diagonalize this R-operator.

\vspace*{2mm}\noindent
MSC:  
\end{abstract}

\maketitle

\section{Introduction}

Quantum groups have become an indispensable 
tool in many areas of mathematical physics and mathematics.
In a broad class of quantum theoretical models it has turned out that
finding a relation to a quantum group is the key for obtaining
exact information about the spectrum or the correlation functions.

So far most of the vast amount of work
devoted to quantum group theory and their physical applications
was concerned with quantum groups that can be studied in a 
purely algebraic manner. This is the case e.g. if the relevant
representations are highest weight representations, as is often
assumed.

However, in physical applications to quantum theoretical models the 
choice of a scalar product on the space of states usually determines
the hermiticity relations for the representatives of the 
quantum group generators. In many cases like 
those corresponding to the so-called {\it non-compact} quantum groups 
it turns out that the corresponding unitary representations are 
always infinite-dimensional and generically neither of highest nor 
lowest weight type. 
In order to exploit the information provided by the appearance of such a 
quantum group it is clearly important to have efficient 
mathematical tools for analyzing the corresponding 
representation theory. 

Unfortunately there are comparatively few results about the 
representation theory of non-compact quantum groups. This seems to
be an important obstacle for making further progress in many
quantum integrable models. Nevertheless there is an
interesting example that 
was first studied independently in \cite{PT2},
\cite{F1} and \cite{PT}. 
These references were considering a particular class of 
infinite dimensional unitary representations of $U_q(\fsl(2,\BR))$. 
This class of representations, henceforth denoted $\CP_s$, $s\in\BR$, 
was the also first to be used in a concrete
physical application: Understanding the tensor products of the 
above-mentioned representations was crucial for obtaining exact results 
on the quantum Liouville theory \cite{PT2}\cite{T}.

A rather remarkable duality phenomenon was observed in \cite{F2},\cite{PT2} 
and \cite{F1}. 
The representations in question are simultaneously 
representations of the two  
quantum groups $U_q(\fsl_2)$ and $U_{\tilde{q}}(\fsl_2)$
with deformation parameters $q=e^{\pi i b^2}$ and
$\tilde{q}=e^{\frac{\pi i}{b^2}}$ respectively. 
This duality turns out to 
be deeply related to the quantum field theoretical self-duality
of Liouville theory under the change of the coupling constant $b$
into $b^{-1}$ \cite{T}. Moreover, it is this duality 
under $b\ra b^{-1}$ that allows one to cover the so-called
{\it strong-coupling} regime where $|b|=1$ by analytic continuation
of the results for real values of $b$ \cite{PT2}\cite{FVK}\cite{T}\cite{FK2}.

The results of the present paper 
clarify the origin of this phenomenon to a certain extent. 
Given the operators $\SX$ 
representing one of the two algebras, say $U_q(\fsl_2)$, 
one may obtain the representatives of the 
second algebra $U_{\tilde{q}}(\fsl_2)$ 
as {\it nonpolynomial} operator functions of the operators $\SX$.
This is found to be consistent with the respective co-products.
In particular, restricting attention to only one of 
the two algebras does not lead to any degeneracy 
as is sometimes suggested in the literature. 

Faddeev has 
proposed to unite the quantum groups $U_q(\fsl_2)$ and 
$U_{\tilde{q}}(\fsl_2)$ into some enlarged object for which he has
coined the name ``modular double''.
The proposal of \cite{F1} as refined
in \cite{KLS} amounts to defining it as the product of 
$\hat{U}_q(\fsl_2)$ and $\hat{U}_{\tilde{q}}(\fsl_2)$,
where, roughly speaking, 
$\hat{U}_q(\fsl_2)$ is obtained from $U_q(\fsl_2)$
by adjoining a sign to the center generated by the Casimir of 
$U_q(\fsl_2)$. We feel that this definition for the modular double
has a disadvantage, though. 
Most representations of $\hat{U}_q(\fsl_2)\ot\hat{U}_{\tilde{q}}(\fsl_2)$ 
are simply tensor products of representations of the two factors.
The representations $\CP_s$ on the contrary are distinguished by the fact 
that they do {\it not} factorize as a tensor product of representations for
$U_q(\fsl_2)$ and $U_{\tilde{q}}(\fsl_2)$. This is what makes the
duality under $b\ra b^{-1}$ a nontrivial statement. Since the 
category formed by the representations $\CP_s$
is {\it closed} under tensor products \cite{PT} it seems 
natural to look for the
group-like object that contains the interesting representations 
$\CP_s$ {\it only}.


We therefore propose to look for a
definition of the ``modular double'' that 
excludes the representations of
$U_q(\fsl_2)\ot U_{\tilde{q}}(\fsl_2)$ which factorize trivially.
This definition should still capture the duality phenomenon
mentioned above. As we have indicated, this naturally leads us to consider 
nonpolynomial functions of the generators.
The basic objects underlying our approach
to the modular double will be an algebra $\CA$ of {\it bounded}
operators, a coproduct $\De$ on $\CA$, an invariant integration
(Haar-measure) on $\CA$ and the R-operator proposed in \cite{F1}. 
We are going to establish the main properties satisfied by these objects,
which are all {\it self-dual} under $b\ra b^{-1}$.

The algebra $\CA$ can be thought of as being generated from
operators that represent $U_q(\fsl_2)$ in a similar 
(in fact, closely related) way as the
algebra of bounded operators on $L^2(\BR)$ is related to the
usual quantum mechanical position and momentum operators $\sx$ and
$\spp$. Our point of view is inspired by the one of 
Woronowicz \cite{W}, which has stimulated the
development of a theory for noncompact quantum groups in a
$C^*$-algebraic framework, see e.g. \cite{KV} and references therein. 
However, although we believe that our results represent substantial
progress towards a proof that the modular double fits into such a 
$C^*$-algebraic framework, it was not our aim to actually carry out 
such a proof here.

We also clarify the relation between 
the R-operator proposed in \cite{F1} and 
the calculus of Clebsch-Gordan and Racah-Wigner-coefficients of
\cite{PT}. Establishing this relation is 
important for the following reason. 
In \cite{T} it was shown that certain families of 
representations of the Virasoro algebra and of the quantum 
group $U_q(\fsl(2,\BR))$ behave equivalently under the
respective product operations (fusion and tensor product).
Together with the results of the present paper it follows that
the respective {\it braiding} operations are equivalent as well.

To be specific we will mostly consider the case that 
the deformation parameter is of the form $q=e^{\pi i b^2}$, where
$b\in(0,1)$. However, our results will carry over to the 
``strong coupling regime'' $|b|=1$, see the remarks in 
Subsection \ref{strong}.

\section{Definitions and main results}
\subsection{Star algebra $U_q(\fsl(2,\BR))$}
$U_q(\fsl(2,\BR))$ is a Hopf-algebra with 
\begin{equation}\label{def}
\begin{aligned}
{}& \text{generators:}\quad  E,\quad F,\quad K,\quad K^{-1};\\
& \text{relations:}\quad 
KE=qEK,\qquad KF=q^{-1}FK,\qquad[E,F]=\frac{K^2-K^{-2}}{q-q^{-1}};\\
& \text{star-structure:}\quad 
K^*=K, \qquad E^*=E, \qquad F^*=F \,.
\end{aligned}\end{equation}
The center of $U_q(\fsl(2,\BR))$ is generated by the $q$-Casimir
\begin{equation}\label{Cas}
C=FE + \frac{qK^2+q^{-1}K^{-2}-2}{(q-q^{-1})^2} \,.
\end{equation}
Compared to the definition used in \cite{PT} we have redefined
$F\ra -F$. This will allow us to realize $F$ by positive
operators. 

\subsection{The representations $\CP_s$ of $U_q(\fsl(2,\BR))$}

In the present paper we will study a one-parameter class 
$\CP_{s}$, $s\in \BR$, of representations of $U_q(sl(2,\BR))$. 
They are constructed as follows: The representation will be 
realized on the space $\CP_{s}$ of entire analytic functions 
$f(x)$ that have a Fourier-transform $\tilde{f}(\om)$ which 
is meromorphic in $\BC$ with the possible poles contained in
\begin{equation}
\CS_s\;\equiv\;
\Big\{\pm\om\,=\,s+i\bigl(\fr{Q}{2}+nb+mb^{-1}\bigr),\;
 n,m\in\BZ^{\geq 0} \,\Big\}, 
\end{equation}
where $Q=b+b^{-1}$.
The representation of $U_q(\fsl(2,\BR))$ on $\CP_s$ is then
defined by choosing the representatives $\pi_s(X)$ for $X=E,F,K$ 
to be the following finite difference operators 
\begin{equation}\label{EFK1}
\begin{aligned}
\pi_{s}(E)\;\equiv\;\SE_{s}\;\equiv\;&e^{+\pi b \sx} \,
 \frac{\cosh\pi b(\spp-s)}{\sin\pi b^2} \, e^{+\pi b \sx}
\\
 \pi_{s}(F)\;\equiv\;\SF_{s}\;\equiv\;&
 e^{-\pi b \sx} \, \frac{\cosh\pi b(\spp+s)}{\sin\pi b^2} \,
 e^{-\pi b \sx}
\end{aligned}
 \qquad \pi_{s}(K)\;\equiv\SK_s\;\equiv\;e^{-\pi b\spp} \,,
\end{equation}
where $\spp$ and $\sx$ are self-adjoint operators satisfying
$[\spp,\sx]=(2\pi i)^{-1}$. By embedding $\CP_s$ as a dense
subspace into the Hilbert space $L^2(\BR)$ one obtains a unitary
representation of $U_q(\fsl(2,\BR))$ generated from the self-adjoint
operators $\SE$, $\SF$ and $\SK$ \cite{S}.

\subsection{The representations $\CH$ and $\CK$ of $U_q(sl(2,\BR))$}
\label{CH}

We will find it convenient to formulate our results in a
``universal'' setting. Let us define 
$\CK\equiv L^{2}(\BR\ti\BR)$. The algebra $\CB(\CK)$ 
of bounded operators on $\CK$ is generated by two pairs 
$(\sx_i,\spp_i)$, $i=1,2$ satisfying 
$[\spp_i,\sx_i]=(2\pi i)^{-1}$. The action of 
$U_q(\fsl(2,\BR))$ on $\CK$ is defined by  
\begin{equation}\label{EFK2}
\begin{aligned}
\pi_{\CK}(E)\;\equiv\;\SE\;=&\;
 e^{\pi b(\sx_\one-\sx_\two)} \,
 \frac{\cosh\pi b \spp_2}{\sin\pi b^2} \,
 e^{\pi b(\sx_\one-\sx_\two)},\\
\pi_{\CK}(F)\;\equiv\;\SF\;=&\;
 e^{\pi b(\sx_\two-\sx_\one)} \, 
 \frac{\cosh\pi b \spp_1}{\sin\pi b^2} \,
 e^{\pi b(\sx_\two-\sx_\one)} \,,
\end{aligned}
\qquad \pi_{\CK}(K)\;\equiv\;
 \SK=e^{\frac{\pi b}{2}(\spp_2-\spp_1)} \,.
\end{equation}
This representation of $U_q(\fsl(2,\BR))$ on $\CK$ is reducible:
The operator $\ss=\frac{1}{2}(\spp_1+\spp_2)$ commutes with 
$\SE$, $\SF$, $\SK$ 
and determines the representation of the Casimir via
\begin{equation}\label{Casimir}
\SC\equiv \pi_{\CK}(C) \,=\, 
\frac{\cosh^2\pi b\ss}{\sin^2\pi b^2} \,. 
\end{equation}
The action of (\ref{EFK2}) on an eigenspace of 
$\ss$ reduces to the action (\ref{EFK1}) on $\CP_s$
upon identification 
$\spp=\frac{1}{2}(\spp_1 - \spp_2)$ and  $\sx=\sx_1 - \sx_2$.
This means that $\CK$ decomposes into the representations 
$\CP_s$ as follows:
\begin{equation}\label{CHdecomp}
\CK\;\simeq\; \int_{\BR}^{\oplus}ds\;\CP_s \,.
\end{equation}
The representations $\CP_s$ and $\CP_{-s}$ are unitarily equivalent \cite{PT}:
There exists a unitary operator $\SJ_s:\CP_s\ra\CP_{-s}$ such that
$\SX_{-s}\SJ_s=\SJ_s\SX_s$ for all $X\in U_q(\fsl(2,\BR))$. The operator
$\SJ_s$ defines an operator $\SJ:\CK\ra\CK$ if one considers
the operator function $\SJ\equiv\SJ_{\ss}$.

It will sometimes
be convenient to consider instead of $\CK$ a space $\CH$ in which
$\CP_s$ and $\CP_{-s}$ are identified:
\[
\CH\;=\;\{ v\in\CK; \,(\id-\SJ)v=0\}.
\]
$\CH$ is of course isomorphic to $\int_{\BR^+}^{\oplus}ds \;\CP_s$.

\subsection{Operator functions of $\SE$, $\SF$, $\SK$}
It is important to also consider nonpolynomial functions of the 
operators $\SE$, $\SF$, $\SK$. Let us first note that standard
functional calculus for positive selfadjoint operators allows one to 
consider complex powers of the generators such as $\SE^{\ga}$, $\ga\in\BC$.
The following result offers a partial explanation for the phenomenon
of modular duality.
\begin{lem}
(i) The operators $\tilde{\SE}_s$, $\tilde{\SF}_s$,  $\tilde{\SK}_s$ obtained
by replacing $b\ra b^{-1}$ in \rf{EFK1},
\begin{equation}\label{EFK'}
\begin{aligned}
\tilde{\SE}_s\;\equiv\;&e^{+\pi b^{-1} \sx} \,
 \frac{\cosh\pi b^{-1}(\spp-s)}{\sin\pi b^{-2}} \, e^{+\pi b^{-1} \sx}
\\
\tilde{\SF}_s\;\equiv\;&
 e^{-\pi b^{-1} \sx} \, \frac{\cosh\pi b^{-1}(\spp+s)}{\sin\pi b^{-2}} \,
 e^{-\pi b^{-1} \sx}
\end{aligned}
 \qquad \tilde{\SK}_s\;\equiv\;e^{-\pi b^{-1}\spp} \,,
\end{equation}
generate a representation of $U_{\tilde{q}}(sl(2,\BR))$ with 
$\tilde{q}=e^{\pi i b^{-2}}$. The generators 
$\tilde{\SE}_s$, $\tilde{\SF}_s$,  $\tilde{\SK}_s$ commute 
with the operators $\SE_s$, $\SF_s$, $\SK_s$ on $\CP_s$.\\
(ii) For $\ga=b^{-2}$ we have
\begin{equation}\label{dualeqs}\begin{aligned}
\big(\sin(\pi b^2)\SE_s\big)^\ga\;=&\;\sin(\pi b^{-2})\,\tilde{\SE}_s,\\
\big(\sin(\pi b^2)\SF_s\big)^\ga\;=&\;\sin(\pi b^{-2})\,\tilde{\SF}_s,
\end{aligned}\qquad
\SK^\ga_s\;=\;\tilde{\SK}_s.
\end{equation}
Being operator functions of $\SE_s$, $\SF_s$, $\SK_s$, the operators 
 $\tilde{\SE}_s$, $\tilde{\SF}_s$,  $\tilde{\SK}_s$ do {\bf not}
commute with $\SE_s$, $\SF_s$, $\SK_s$ in the usual sense 
(commutativity of the spectral projections).
\end{lem}

We shall now define an algebra of bounded operators that can be considered
as operator functions of $\SE$, $\SF$, $\SK$.  
To begin with, let 
$\CO_s$, $s\in\BR$ be a family of bounded operators on $L^2(\BR)$
such that 
\[ \sup_{s\in\BR}\;\lVert \CO_s \rVert < \infty.
\]
A bounded operator $\CO$ on $\CK$ can be defined for each such family 
$(\CO_s)_{s\in\BR}$ by means of \rf{CHdecomp}. 
These operators $\CO$ form a subalgebra $\CB_0$ of the algebra of 
all bounded operators on $\CH$. Let $\CB$ be 
the $C^*$ subalgebra obtained as the completion of $\CB_0$ 
w.r.t. the operator norm. This algebra can be thought of as being generated
from the unbounded elements $\spp$, $\sx$, $\ss$.

However, there is no canonical way to define ${\rm sgn}(\ss)$ as a function 
of $\SE$, $\SF$, $\SK$. The center of the algebra of bounded operators 
generated from $\SE$, $\SF$, $\SK$ should be generated from 
operator functions of the Casimir, or equivalently $|\ss|$, cf. eqn.
\rf{Casimir}. This is closely 
related to the fact that the representations $\CP_{s}$ and $\CP_{-s}$
are unitarily equivalent. Elements of the ``true''
algebra $\CA\subset\CB$ should therefore commute with the operator 
$\SJ$ which establishes the equivalence between $\CP_{s}$ and $\CP_{-s}$,
\begin{equation}
\CA\;\equiv\;\{\SO\in\CB; \SJ^{-1}\SO\SJ=\SO\}.
\end{equation}
This amounts to considering only those elements
of $\CB$ that leave $\CH$ invariant.

\subsection{The Hopf algebra structure}

A co-product is defined on $U_q(\fsl(2,\BR))$ via
\begin{equation}\label{De}
\begin{aligned}\De(E)=&E\ot K+K^{-1}\ot E \,,\\
\De(F)=&F\ot K+K^{-1}\ot F \,,
\end{aligned}
 \qquad \De(K)=K\ot K \,. 
\end{equation}
In the following we shall adopt the convention to denote 
\begin{equation}\label{conv}
 \De(\SX)\;\equiv\;(\pi_{\CH}\ot\pi_{\CH})\circ
 \De(X)\quad\text{for}\quad X\in \CU_q(\fsl(2,\BR)) \,.
\end{equation}
It follows from  \cite[Theorem 2]{PT} that $\De(\SE)$,
$\De(\SF)$ and $\De(\SK)$ are self-adjoint and positive and
therefore generate a representation of $U_q(\fsl(2,\BR))$ 
on $\CH\ot \CH$. The following Lemma proven in Section 3
establishes consistency of the co-product with 
modular duality:
\begin{lem}\label{codual}
The co-product (\ref{De})\rf{conv} acts on the 
dual part of the modular double as follows.
\begin{equation}\label{ded}\begin{aligned}
\De(\tilde{\SE})\;&=\;\tilde{\SE}\ot\tilde{\SK}+\tilde{\SK}^{-1}\ot\tilde{\SE},
\\
\De(\tilde{\SF})\;&=\;\tilde{\SF}\ot\tilde{\SK}+\tilde{\SK}^{-1}\ot\tilde{\SF},
\end{aligned}\qquad 
\De(\tilde{\SK})\;=\;\tilde{\SK}\ot\tilde{\SK}.
\end{equation}
\end{lem}

A representation of the co-product on the algebra $\CA$ can be defined
by means of the Clebsch-Gordan maps defined in \cite{PT}. 
These maps yield a three parameter family of maps 
$\SC[s_3|s_2,s_1]:\CP_{s_\two}\ot\CP_{s_\one}\ra \CP_{s_\three}$
that satisfy the intertwining property
\begin{equation} \label{intertw}
\SC[s_3|s_2,s_1]\;\circ
(\pi_{s_{\two}}\ot \pi_{s_{\one}})\circ\De(X)\,=\,
\pi_{s_\three}(X)\circ\SC[s_3|s_2,s_1]\;
  \,
\end{equation}
and extends to a two-parameter family of unitary operators
$\SC[s_2,s_1]:L^{2}(\BR^2)\;\ra\; 
\CH$. 
Let us introduce the
operators $\ss_\one=\id\ot\ss$, $\ss_{\two}=\ss\ot\id$ on $\CH\ot\CH$
respectively. The identification \rf{CHdecomp} allows us to 
consider $\SC[\ss_{\two},\ss_\one]$ as a unitary
operator 
\[ 
\SC\;\;:\;\; 
\CH\ot\CH\ra \CH \ot \CH_{\rm spec}^{(2)}\ot
\CH_{\rm spec}^{(1)},
\]
where the operators $\ss_i$, $i=1,2$  are realized
on the spaces $\CH_{\rm spec}^{(i)}\simeq L^2(\BR^+)$
as multiplication operators.
 
For each element $\SX\in\CA$ we may now define $\De(\SX)$ by
\begin{equation}
\De(\SX)\;\equiv\;\SC^{\dagger}\circ \big(\SX\ot \id\ot\id\big)\circ \SC.
\end{equation}
Since $\SC$ is unitary and $\SX$ is bounded we clearly have boundedness
of $\De(\SX):\CH\ot\CH\ra\CH\ot\CH$.
\begin{thm}\label{coass}
The coproduct $\De$ is coassociative on $\CA$, i.e.
\[
(\id\ot\De)\circ\De(\SX) \;=\; {(\De\ot\id)}\circ{\De(\SX)}\quad
\text{for any}\quad \SX\in\CA. 
\]
\end{thm}

The antipode consistent with (\ref{De}) is defined as an 
anti-automorphism of $U_q(\fsl(2,\BR))$ such that
\begin{equation}
\sigma(K)= K^{-1} \,,\qquad 
\sigma(E)= - q E \,, \qquad
\sigma(F)= -q^{-1} F \,.
\end{equation}
The action of the antipode on nonpolynomial functions of $\SE$, $\SF$ 
and $\SK$ can be introduced by means of
\begin{equation}
\sigma(\spp)\;=\;-\spp,\qquad \sigma(\ss)\;=\;-\ss,\qquad
\sigma(\sx)\;=\;\sx+\fr{i}{2}Q.
\end{equation}
The fact that $\sx$ is shifted by an imaginary amount means that 
$\sigma$ is not defined on all of $\CA$. This unboundedness
of the antipode is not unexpected \cite{KV}.

\subsection{The Haar-measure}

Let us first note that the decomposition \rf{CHdecomp} induces
a family of projections $\pi_s:\CA\ra \CB(L^2(\BR))$. We shall often
use the shorthand 
notation $\SO_s\equiv \pi_s(\SO)$.
\begin{defn}
Define linear functionals $\sh_{\rm l}$ and $\sh_{\rm r}$ on 
dense subsets 
$\CA_{\sh}^{\rm l}$ and $\CA_{\sh}^{\rm r}$ of $\CA$ respectively by
\begin{equation}\label{Hm}
\begin{aligned}
\sh_{\rm l}(\SO)\;=\;&
\int_{0}^{\infty}d{\rm m}(s)\;{\rm Tr}(e^{-2\pi Q\spp}\SO_s),\\ 
\sh_{\rm r}(\SO)\;=\;&\int_{0}^{\infty}
d{\rm m}(s)\;{\rm Tr}(e^{+2\pi Q\spp}\SO_s), 
\end{aligned}\end{equation}
where the measure m is defined by 
\begin{equation}
d{\rm m}(s)\;\equiv\;4\sinh2\pi bs\sinh2\pi b^{-1}s \;ds.
\end{equation}
\end{defn}
\begin{thm} \label{invthm}
(i) The Haar-measures $\sh_{\rm l}$ and $\sh_{\rm r}$ 
are left and right invariant respectively,
\begin{equation}
\begin{aligned}\label{Hinv1}
(\id\ot\sh_{\rm l})\circ\De(\SO)\;=\;\sh_{\rm l}(\SO)\,\id\;,\\
(\sh_{\rm r}\ot\id)\circ\De(\SO)\;=\;\sh_{\rm r}(\SO)\,\id\;,
\end{aligned}
\end{equation}
where we assume $\SO$ to be taken from the respective domains 
of definition.\\
(ii) For any $\SX \in U_q(\fsl(2,\BR))$, the Haar-measures satisfy,
respectively
\begin{equation}\label{Hinv2}
\begin{aligned}
 \sh_{\rm l} ( {\rm ad^l_\SX} \SO) \;=\;
 \sh_{\rm l}(\SO)\ \epsilon(\SX) \,,\\
 \sh_{\rm r} ( {\rm ad^r_\SX} \SO) \;=\;
 \sh_{\rm r}(\SO)\ \epsilon(\SX) \,,
\end{aligned}
\end{equation}
where $\epsilon(\SX)$ is the co-unit, and
the left and right $q$-adjoint actions are defined as
${\rm ad^l_X} (Y) = \sum_i X^\prime_i Y \sigma(X^{\prime\prime}_i)$
and 
${\rm ad^r_X} (Y) = \sum_i \sigma(X^\prime_i) Y X^{\prime\prime}_i$
if $\Delta(X)=\sum_i X^\prime_i \otimes X^{\prime\prime}_i$.
\end{thm}
\begin{rem}
We believe that the triple $(\CA,\De,h)$ that we have defined above
constitutes a somewhat more satisfactory definition of the modular double,
although more work is needed to show that it fits into the 
axiomatics for noncompact quantum groups of \cite{KV}. The self-duality
under $b\ra b^{-1}$ is manifest in this formulation.

It also becomes clear that the modular double can not be considered
as a deformation of a classical group: The Haar-measure has 
no classical limit $b\ra 0$ due to the factor $Q=b+b^{-1}$ 
that appears in the definition of $\sh$.
\end{rem}

\subsection{The R-operator}

To begin with, we introduce the special function $g_b(x)$ that will be 
used to define the R-operator. It may be defined via (recall that
$Q=b+b^{-1}$)
\begin{equation}\label{gG}
\log g_b(x)\;=\;-\int\limits_{\BR+i0}\frac{dt}{t}\;\frac{e^{^\frac{tQ}{2}}\,
x^{^\frac{t}{2\pi i b}}}{(1-e^{bt})(1-e^{t/b})}.
\end{equation}
Let us furthermore introduce an anti-self-adjoint element $\SH$ such
that $K=q^H$. Define
\begin{equation}\label{Rdef}
 \SR\,=\, q^{\SH\ot\SH} \,
 g_b\bigl(4 (\sin\pi b^2)^2 \SE\ot
 \SF\bigr)\, q^{\SH\ot\SH} \,,
\end{equation}
where $\SH \equiv \pi_{\CH}(H)$.
As we will explain below (see Corollary~\ref{Rprorep}), 
$\SR$ coincides with the R-operator proposed 
by L.~Faddeev in~\cite{F1}. 
Notice that the property (\ref{Gcc}) implies that
$|g_b(x)|=1$ for $x\in{\BR}^+$. This makes
$\SR$ manifestly unitary.

\begin{thm}\label{T1}
The operator $\SR$ has the following properties:
\begin{eqnarray}
 & {\rm (i)} &  \SR\,\De(\SX) \,=\, 
     \De'(\SX)\,\SR \,,  \label{i} \\
 & {\rm (ii)} & (\id\ot\De)\SR\,=\, \SR_{13}\SR_{12} ,\qquad
 (\De\ot\id)\SR\,=\, \SR_{13}\SR_{23} \,, \label{ii} \\
 & {\rm (iii)} & (\sigma \ot \id)\SR\,=\, \SR^{-1} , \qquad 
 (\id \ot \sigma)\SR\,=\, \SR^{-1} , \qquad
 (\sigma \ot \sigma)\SR\,=\, \SR \,.  \label{iii}
\end{eqnarray}
\end{thm}

\begin{rem}
The R-operator allows us to introduce the braiding of tensor 
products of the representations $\CP_{s}$. Specifically, let the 
operator $\SB:\CP_{s_{\two}}\ot \CP_{s_\one}\ra \CP_{s_{\one}}\ot 
\CP_{s_\two}$
be defined by $\SB_{s_{\two},s_\one}
\equiv \SP\SR_{s_{\two},s_\one}$, where $\SP$ is the 
operator that permutes the two tensor factors. Property 
(i) {}from Theorem~\ref{T1} implies as usual that 
$\SB_{s_{\two},s_\one}\circ \De(\SX)=
 \De(\SX)\circ\SB_{s_{\two},s_\one}$.
\end{rem}

\subsection{Integral operator representation}

The operator $\SR$ can clearly be projected to an operator
$\SR_{s_\two s_\one} \equiv
 (\pi_{s_\two} \otimes \pi_{s_\one}) \, \SR$ 
on $\CP_{s_{\two}}\ot \CP_{s_{\one}}$.
The action of this operator admits
a representation by means of a distributional kernel: 
\begin{thm}\label{Tint}
Let $\tilde{\psi} (k_2,k_1) =  \int_\BR dx_2 dx_1 
 e^{2\pi i (k_1 x_1 + k_2 x_2)}\psi (x_2,x_1)$ be a 
Fourier transform of 
$\psi (x_2,x_1) \in \CP_{s_{\two}}\ot \CP_{s_{\one}}$.
The action of the $\SR$-operator on 
\hbox{$\CP_{s_{\two}}\ot \CP_{s_\one}$} admits
the following representations as an integral operator in
``coordinate'' and ``momentum'' space respectively:
\begin{eqnarray}
 && \label{io1}
 \bigl( \SR_{s_2s_1} \, \psi \bigr)(x_2,x_1) \,=\, 
 \int_\BR dx'_2 dx'_1 \, 
 \SR_{s_2s_1}(x_2,x_1|x'_2,x'_1) \, \psi(x'_2,x'_1) \,, \\
 && \label{io2} 
 \bigl( \SR_{s_2s_1} \, \tilde{\psi} \bigr)(k_2,k_1) \,=\, 
 \int_\BR dk'_2dk'_1 \, 
 \tilde{\SR}_{s_2s_1}(k_2,k_1|k'_2,k'_1) \,  
 \tilde{\psi}(k'_2,k'_1) \,,
\end{eqnarray}
with the kernels given by
\begin{equation}\label{ker1}
\begin{aligned}
  \SR_{s_2s_1} &(x_2,x_1 | x'_2,x'_1) = 
 e^{2\pi i\left( s_1(x'_1-x_1) + s_2(x_2-x'_2) 
 + \frac{i Q}{2} (x_2 + x'_2 - x_1 - x'_1) + s_1 s_2
 + \frac 14 Q^2 \right)}  \\
 & \times
 \frac{G_b\bigl(\frac Q2 +\frac i2 (s_1+s_2)+ i(x_2-x_1)\bigr)}%
 {G_b\bigr(Q+ \frac i2 (s_1-s_2) + i(x_2 - x'_1)\bigr)} \,
 \frac{G_b\bigl(\frac Q2 -\frac i2 (s_1+s_2) + i(x'_2-x'_1)\bigr)}%
 {G_b\bigl(Q +\frac i2 (s_2-s_1) + i(x'_2-x_1)\bigr)} \,,
\end{aligned}
\end{equation}
\begin{equation}\begin{aligned}\label{ker2}
 \tilde{\SR}_{s_2s_1} & (k_2,k_1| k'_2,k'_1) =\\ 
 &=\delta(k'_2 + k'_1 - k_2 - k_1) \,
 \frac{e^{-\pi i (k'_1 k_2 + k_1 k'_2)}}{G_b(Q+i(k'_1 - k_1))}\, %
 \frac{w_b(s_1 + k_1)}{w_b(s_1 + k'_1)} \,
 \frac{w_b(s_2 - k_2)}{w_b(s_2 - k'_2)} .
\end{aligned}
\end{equation}
The functions $w_b(x)$ and $G_b(x)$ are close relatives of the
function $g_b(x)$ that will be defined in Subsection \ref{specfkt} below,
and $G^{-1}(Q+ix)$ is taken as a short notation for the 
distribution $G^{-1}\big(Q+i(x+i0)\big)$.
\end{thm}

\subsection{Highest weight representations}
\label{FDR}

In order to demonstrate that the R-operator we are considering
here indeed deserves to be called ``universal'' we are now going
to show that the usual R-matrix for highest weight representations
of $U_q(\fsl(2,\BR))$ can be extracted from the analytic properties
of the matrix elements given in Theorem \ref{Tint}.

As a preparation let us consider the representation of $U_q(\fsl(2,\BR))$ on
the dual space $\CP'_s$ of $\CP_s$. An interesting class of 
elements of $\CP'_s$ is furnished by the (complexified) 
delta-functionals $\de_{k}$,
\[ \bra\de_k,f\ket=f(k).
\]
The $\de_{k}$ are well-defined for all $k\in\BC\setminus\CS_s$, and the  
action of $U_q(\fsl(2,\BR))$ is realized by
\begin{equation}\label{deltatrf}
\begin{aligned}
\SE_s^t\de_{k}\;=\;&+\big[\fr{Q}{2b}-\fr{i}{b}(k-s)\big]_q\de_{k+ib}\\
\SF_s^t\de_{k}\;=\;&-\big[\fr{Q}{2b}+\fr{i}{b}(k+s)\big]_q\de_{k-ib}
\end{aligned}\qquad \SK_s^t\de_{k}\;=\;e^{-\pi bk}\de_{k},
\end{equation}
where $[t]_q \equiv \frac{\sin(\pi b^2 t)}{\sin(\pi b^2)}$ is the 
standard definition of a q-number, and the superscript ``$t$''
on the generators indicates transposition.

Let us restrict attention to the set $\CD_s$ of functionals $\de_k$ for 
which $k$ is an element of 
$\big\{ k=-s+i\big(\frac{Q}{2}+nb\big), n\in\BZ^{\geq 0}\big\}$. 
It is easy to verify that \rf{deltatrf} realizes a highest
weight representation on $\CD_s$.

\begin{thm}\label{T2} The 
action of $\SR^t$ 
on $\CH\ot\CD_{s}$ is given by
\begin{equation}\label{Rkk}
\begin{aligned}
 \SR_s^+ =& \, q^{\SH\ot\SH_s} \,
 \sum_{n=0}^{\infty} \frac{q^{\frac{1}{2}(n^2-n)}}%
 {\prod_{k=1}^n [k]_q} \, \bigl((q-q^{-1})\SE\ot\SF_s\bigr)^n \, 
 q^{\SH\ot\SH_s} \,.
\end{aligned}
\end{equation}
\end{thm}

\subsection{Diagonalization of the R-operator}

\begin{thm}\label{T3} The Clebsch-Gordan maps 
$\SC[s_3|s_2,s_1]$ diagonalize the R-operator in the 
following sense: 
\begin{equation}\label{Rdiag}
\SC[s_3|s_1,s_2] \,\SB_{s_{\two}s_{\one}}^{}\,=\,
\Om(s_3|s_2,s_1) \, \SC[s_3|s_2,s_1] \,, 
\end{equation}
with eigenvalue $\Om(s_3|s_2,s_1)$ given as
\[
\Om(s_3|s_2,s_1)\;=\; 
e^{-\pi i(h_{s_\three}-h_{s_\two}-h_{s_\one})},\qquad
h_s\;\equiv \;s^2-\fr{Q^2}{4}.\]
\end{thm}

\subsection{The strong coupling regime $|b|=1$} \label{strong}

We would finally like to point out that our results carry over to
the strong coupling regime $b=e^{i\theta}, \theta\in [0,\pi/2)$. 
This is almost
obvious for those results whose proof relies mainly on the
properties of the special functions $g_b(x)$, $G_b(x)$ and $w_b(x)$.
In this case the operators $\SE$, $\SF$ and $\SK$ are {\it normal}
(as follows from eqn. \rf{EFw} below), and the hermitian conjugation
acts as 
\[ 
\SE^{\dagger}\;=\;\tilde{\SE},\qquad
\SF^{\dagger}\;=\;\tilde{\SF},\qquad
\SK^{\dagger}\;=\;\tilde{\SK}.
\]

Concerning the results that rely on \cite{PT} 
one may note that they all amount to certain identities between
distributions that are defined by a standard analytic regularization
in terms of the meromorphic functions $g_b(x)$, $G_b(x)$ and $w_b(x)$.
The relevant analytic properties underlying the
validity of these identities all remain intact upon
analytically continuing from the case of real $b$ to $|b|=1$.

\section{Preliminaries and auxiliary results}
\setcounter{equation}{0}

\subsection{Special functions}\label{specfkt}

The Double Gamma function $\Gamma_2(x|\omega_1,\omega_2)$
was introduced and studied by Barnes \cite{Ba}. 
In what follows we will be dealing with (recall that
$Q=b+b^{-1}$) 
\begin{equation}\label{Gb}
 G_b(x) \equiv e^{ \frac{\pi i}{2} x(x-Q)} \,
 \frac{\Gamma_2(x|b^{-1},b)}{\Gamma_2(Q-x|b^{-1},b)} \,.
\end{equation}
This function is closely related to 
the remarkable special functions introduced under the
names of ``quantum dilogarithm'' in \cite{FK1} and 
``quantum exponential function'' in \cite{W2}.
$G_b(x)$ is a meromorphic function that has the 
following properties \cite{Ba,Sh}:
\begin{eqnarray}
\text{self-duality} \label{Gdual}  && 
 G_b(x)=G_{b^{-1}}(x) \,,    \\
\text{functional equation}  \label{Gfunrel} &&
 G_b(x+b)=(1-e^{2\pi i bx}) \, G_b(x) \,,  \\
\text{reflection property} \label{Grefl} &&
 G_b(x) \, G_b(Q-x)=e^{\pi i x(x-Q)} \,,  \\
   \label{Gcc} \text{complex conjugation} &&
  \overline{G_b(x)} = 
 e^{\pi i \bar{x}(Q - \bar{x})} \, G_b(\bar{x}) \,, \\
\text{asymptotics}  \label{Gasym} &&
  G_b(x) \sim
\left\{
\begin{array}{ll} 
{}  \overline{\zeta}_b &
\text{ for }\ \Im(x)\ra +\infty \\ [0.5mm]
    \zeta_b \, e^{\pi i x(x-Q)} &
\text{ for }\ \Im(x)\ra -\infty
\end{array}
\right.  \,, \\
 \label{Gan} 
 \begin{array}{r} 
 G_b(x) \ \text{has poles at} \\ [-0.25mm]
  G_b(x) \ \text{has zeros at}
\end{array}
  && 
 \begin{array}{l} 
 x = -nb-mb^{-1} \\ [-0.25mm]
 x = Q+nb+mb^{-1}
\end{array}
 \ \  n,m \in {\BZ}^{\geq 0} \,,
\end{eqnarray}
where $\zeta_b = e^{\frac{\pi i}{4} 
 +\frac{\pi i}{12} (b^2+b^{-2})}$.
By Proposition 5 in \cite{Sh}, the $G_b$-function 
admits for $\Im b^2>0$
the following infinite product representation
\begin{equation}\label{Gprod}
 G_b(x) = \overline{\zeta}_b \ \frac{\prod_{n=1}^\infty 
  (1 - e^{2 \pi i b^{-1} (x -nb^{-1})} ) }%
{\prod_{n=0}^\infty (1 - e^{2 \pi i b (x +nb)} )} \,.
\end{equation}

We are also going to use two other functions that are closely related to
$G_b(x)$, namely
\begin{equation}\label{wb}
g_b(x)\equiv \frac{\overline{\zeta}_b}{G_b\bigl(\frac{Q}{2}+
  \frac{1}{2\pi ib}\log x\bigr)}, \quad\text{and}\quad
w_b(x) \equiv e^{ \frac{\pi i}{2} 
   \bigl( \frac{Q^2}{4} +x^2 \bigr)} \, 
 G_b \Bigl(\frac{Q}{2} - i x \Bigr) \,.
\end{equation} 
The representation \rf{gG} for $g_b$ introduced above follows easily 
from the integral representation for the Double Gamma function introduced
in \cite{Sh}.

For the reader's convenience we shall also list the relevant properties of 
$w_b(x)$ that follow from (\ref{Gdual})--(\ref{Gan}).  
\begin{eqnarray}
\text{self-duality} \label{wdual}  && 
 w_b(x)=w_{b^{-1}}(x) \,,  \\
\text{functional equation}  \label{wfunrel} &&
 w_b(x + i b)= 2 \, w_b(x) \, 
 \sin\pi b \bigl(\fr{Q}{2} - i x\bigr) \,, \\
\text{reflection property} \label{wrefl} &&
 w_b(x) \; w_b(-x) = 1  \,, \\ [0.5mm]
  \label{wcc} \text{complex conjugation} &&
  \overline{w_b(x)} = w_b(-\bar{x}) \,,\\
 \label{wan} 
 \begin{array}{r} 
 w_b(x) \ \text{has poles at} \\ [-0.25mm]
  w_b(x) \ \text{has zeros at}
\end{array}
  && 
 \begin{array}{l} 
 x = - i \, (\frac{Q}{2} +nb+mb^{-1}) \\ [0.25mm]
 x = i \, (\frac{Q}{2} +nb+mb^{-1})
\end{array}
 \ \  n,m \in {\BZ}^{\geq 0} \,.
\end{eqnarray}
Notice that $|w_b(x)|=1$ if $x$ is real. Hence $w_b(\SX)$ is
unitary if $\SX$ is a self-adjoint operator.

\subsection{Operator algebraic preliminaries}

\begin{lem}\label{uvg}
Let $A$ and $B$ be self-adjoint operators such that
$[A,B] = 2 \pi i$. Let $\varphi(t)$ be a
function on the positive real axis and let $\gamma = \frac{1}{b^2}$.
Then we have
\begin{eqnarray}
\label{phiAB}
 & \varphi (u+v) \,=\, w_b\bigl(2\pi(A-B) \bigr) \,
 \varphi \bigl( e^{\frac{b}{2} (A+B)} \bigr) \,
 w_b\bigl(2\pi(B-A) \bigr) \,, & \\
\label{qF}
 & (u + v)^\gamma  \,=\, u^\gamma + v^\gamma \,, &
\end{eqnarray}
where $u = e^{b A}$, $v = e^{b B}$.
\end{lem}

\begin{proof}
It is convenient to introduce $p \equiv \frac{B-A}{2\pi}$ 
and $x \equiv \frac{A+B}{4\pi}$. Observe that
$[p,x]= \frac{1}{2\pi i}$; so that we have
\begin{equation}\label{fe}
 f\bigl( p \bigr) \, e^{\pi b x} =   
 e^{\pi b x} \, f\bigl(p - i \fr{b}{2}\bigr) \,.
\end{equation}
for any function $f(t)$ that is bounded and analytic in the
strip $\frac{b}{2}\leq\Im(p)\leq 0$. 
Using the Baker-Campbell-Hausdorff formula
for Weyl-type operators and the properties of the
function $w_b$, we may calculate as follows:
\[
\begin{aligned}
w_b(-p) \, e^{2\pi b x} \, w_b(p) & \={wrefl}  
\frac{1}{w_b(p)} \, e^{2\pi b x} \, w_b(p) \={fe}
e^{\pi b x} \, 
 \frac{w_b(p+ i\frac{b}{2})}{w_b(p - i\frac{b}{2})}e^{\pi bx}\\
& \={wfunrel}  2 e^{\pi b x} \, (\cosh {\pi b p}) \, e^{\pi b x} 
= e^{\frac b4 (A+B)} \bigl( e^{\frac b2 (B-A)} +
 e^{\frac b2 (A-B)} \bigr) e^{\frac b4 (A+B)}\\
& = u+v \end{aligned}
\] 
The last expression is therefore unitarily equivalent to 
the positive self-adjoint operator $e^{2\pi b x}$. 
Our claim (\ref{phiAB}) follows by applying the standard functional 
calculus of self-adjoint operators.

Relation (\ref{qF}) can be proven along the same lines
taking into account that, thanks to self-duality (\ref{wdual}), 
$w_b$ obeys also the equation 
\begin{equation}\label{wdf}
w_b(x + \fr{i}{b})= 2 \, w_b(x) \, 
 \sin\bigl( \fr{\pi}{b} \bigl(\fr{Q}{2} - i x\bigr) \bigr)\,.
\end{equation}
Therefore, for $\varphi(t)=t^\gamma$ we have
\begin{equation}\nonumber
\begin{aligned}
 (u+v)^\gamma &= w_b(-p) \, e^{\frac{2\pi}{b} x} \, w_b(p)
    \={fe} e^{ \frac{\pi}{b} x} \, 
 \frac{w_b(p+\frac{i}{2b})}{w_b(p-\frac{i}{2b})} e^{ \frac{\pi}{b} x}\\
 & \={wdf}  2 e^{\frac{\pi}{b} x} \, (\cosh {\fr{\pi}{b} p}) \, 
 e^{\frac{\pi}{b} x}
 =  e^{\frac{1}{4b} (A+B)} \bigl( e^{\frac{1}{2b} (B-A)} +
 e^{\frac{1}{2b} (A-B)} \bigr) e^{\frac{1}{4b} (A+B)}
 = u^\gamma + v^\gamma \,.
\end{aligned}
\end{equation}
\end{proof}

\begin{rem}
Another way to prove relation (\ref{qF}) in Lemma~\ref{uvg} 
is to use the $b$-binomial
formula~(\ref{bbin}) that we derive in Appendix. When
$t$ approaches the value $-i \gamma$, the $b$-binomial 
coefficient (\ref{bbc}) vanishes unless $\tau$ takes
special values determined by~(\ref{Gan}). Furthermore,
for $t= -i \gamma$ the $b$-binomial coefficient has 
nonvanishing residues only at $\tau=0$ and $\tau=-i \gamma$.  
The contributions {}from these two poles yield the two terms 
on the r.h.s.~of~(\ref{qF}). 
Similar consideration for $t$ approaching $-i n \gamma$,
$n>1$ shows that the $b$-binomial coefficient has 
nonvanishing residues at $\tau=0, -i\gamma, \ldots, -in\gamma$.
Therefore \hbox{$(u+v)^{n \gamma}$} can be represented as 
sum of $(n+1)$ terms which is analogous to the 
\hbox{$q$-binomial} formula in the compact case.
\end{rem}

The proven Lemma leads to useful representations for the generators and
the R-operator of~$U_q(sl(2,\BR))$. 
For brevity, we denote $\se_b \equiv (2 \sin\pi b^2) \SE$ 
and $\sf_b \equiv (2 \sin\pi b^2) \SF$, whereas 
$\se_{\frac 1b}$ and $\sf_{\frac 1b}$ will stand for
their counterparts with $b$ replaced by $\frac 1b$.
\begin{lem}\label{Rprod}
$\pi_{\CH}(E)$ and $\pi_{\CH}(F)$ 
admit the following representation:
\begin{equation}\label{EFw}
 \se_b = w_b(\spp_2) \, e^{2\pi b(\sx_\one-\sx_\two)} \, 
 w_b(-\spp_2) \,, \qquad
 \sf_b = w_b(\spp_1) \, 
 e^{2\pi b(\sx_\two-\sx_\one)} \,  w_b(-\spp_1) \,.
\end{equation}
$\SR$ may be represented as follows:
\begin{equation}\label{Rpr}
\begin{aligned}
\SR\,=\, q^{\SH\ot\SH} \, 
 \bigl( w_b(\spp_2)\ot w_b(\spp_1) \bigr) \, g_b & 
\bigl( e^{2\pi b(\sx_1-\sx_2)}  \ot 
 e^{2\pi b(\sx_2-\sx_1)} \bigr)\,\\
& \qquad\quad\cdot 
 \bigl( w_b(-\spp_2)\ot w_b(-\spp_1) \bigr)\,
 q^{\SH\ot\SH} \,.
\end{aligned} \end{equation}
\end{lem}
\begin{proof}
In Lemma~\ref{uvg}, we can identify
$\frac{1}{2\pi} (A-B) = \spp_{n+1}$ and 
$\frac{1}{4\pi} (A+B) = \sx_n-\sx_{n+1}$,
where $n=1,2$ (with convention that $n+2\equiv n$).
Then, as seen from the definition (\ref{EFK2}), we have 
$u+v= \se_b$ for $n=1$ and $u+v=\sf_b$ for $n=2$.
Therefore, (\ref{EFw}) is just a particular case
of~(\ref{phiAB}). Furthermore, for functions $\varphi(t)$ defined on $\BR^+$
we have 
\begin{equation}\label{EFphi}
 \varphi(\se_b) = w_b(\spp_2) \, \varphi\bigl(
  e^{2\pi b(\sx_\one-\sx_\two)} \bigr) \, 
 w_b(-\spp_2) \,, \quad
 \varphi(\sf_b) = w_b(\spp_1) \, \varphi\bigl(
 e^{2\pi b(\sx_\two-\sx_\one)} \bigr) \, w_b(-\spp_1) \,.
\end{equation} 
In particular, we can take $\varphi(t)=g_b(x)$
(recall that $|g_b(x)|=1$ for $x\in{\BR}^+$).
For this choice of $\varphi(t)$, the representation 
(\ref{Rpr}) for $\SR$ follows immediately from the 
definition~(\ref{Rdef}). 
\end{proof}

\begin{cor}\label{efbb}
For $\gamma = \frac{1}{b^2}$ we have
\begin{equation}\label{efdual}
 (\se_b)^{\gamma} \,=\, \se_{\frac{1}{b}} \,, \qquad
 (\sf_b)^{\gamma} \,=\, \sf_{\frac{1}{b}} \,.
\end{equation}
\end{cor}
\begin{proof} Notice that, if $u$ and $v$ are identified 
as in the proof of Lemma~\ref{Rprod}, then
\hbox{$u^{\gamma}+v^{\gamma}= \se_{\frac{1}{b}}$} for $n=1$ and 
$u^{\gamma}+v^{\gamma}=\sf_{\frac{1}{b}}$ for $n=2$.
Thus, relations (\ref{efdual}) are a particular
case of~(\ref{qF}) .
\end{proof}

This proves the relations \rf{dualeqs} from Lemma 1. 

\begin{cor}\label{Rprorep}
The definition of the $\SR$-operator proposed in \cite{F1}, 
\begin{equation}\label{RF}
 \SR = q^{\SH \ot \SH} \ \frac{\prod_{n=0}^\infty 
  (1 + q^{2n+1} \, \se_b \otimes \sf_b ) }%
 {\prod_{n=0}^\infty (1 + \tilde{q}^{2n+1} \, 
 \se_{\frac{1}{b}} \otimes \sf_{\frac{1}{b}} )} 
 \ q^{\SH \ot \SH} \,,\quad\tilde{q}=e^{-i\pi b^{-2}},
\end{equation}
which is valid for $b=e^{i\vartheta}$,
$\vartheta\in (0,\frac{\pi}{2})$, coincides
with our definition (\ref{Rdef}).
\end{cor}
\begin{proof}
\[
\begin{aligned}
 {} & q^{-\SH\ot\SH}  \, \SR  \ q^{-\SH\ot\SH}  \={Rdef} 
 g_b\bigl( \se_b \ot \sf_b \bigr)
 \={gG} \overline{\zeta}_b \,  
 \Bigl( G_b\bigl(\fr{Q}{2} +
  \fr{1}{2\pi ib}\log(\se_b \ot \sf_b) \bigr) 
 \Bigr)^{-1} \\
 & \phantom{=} \={Gprod} \frac{\prod_{n=0}^\infty 
 (1 - e^{2 \pi i b (\frac{Q}{2} +nb)} \, \se_b \ot \sf_b)}
 {\prod_{n=1}^\infty   (1 - e^{2 \pi i b^{-1} 
   (\frac{Q}{2} -nb^{-1})} \,(\se_b \ot \sf_b)^{\frac{1}{b^2}} ) }
 \={efdual} \frac{\prod_{n=0}^\infty 
  (1 + q^{2n+1} \, \se_b \otimes \sf_b ) }%
 {\prod_{n=0}^\infty (1 + \tilde{q}^{2n+1} \, 
 \se_{\frac{1}{b}} \otimes \sf_{\frac{1}{b}} )} \,.
\end{aligned}
\]
\end{proof}

\begin{cor}
$\pi_{\CH}(E)$ and $\pi_{\CH}(F)$ 
admit the following representation:
\begin{equation}\label{expEF}
 \se_b = 
 e^{2\pi b(\sx_\one-\sx_\two) +  b \psi_b(\spp_2)} \,, \qquad
 \sf_b = 
 e^{2\pi b(\sx_\two-\sx_\one) +  b \psi_b(\spp_1)}  \,,
\end{equation}
where $\psi_b(t)\equiv i  \, \partial_t (\log w_b(t))$.
\end{cor}
\begin{proof} 
Eqs.~(\ref{EFphi}) for $\varphi(t)=\log (t)$ yield
\begin{equation}\nonumber
\begin{aligned}
 \log(\se_b )  {} &= w_b(\spp_2) \, 
  2 \pi b(\sx_\one-\sx_\two) \, w_b(-\spp_2) 
 = 2\pi b(\sx_\one-\sx_\two) \\ 
 {} & + {} 2\pi b \, [w_b(\spp_2),(\sx_\one-\sx_\two)] \,
 w_b(-\spp_2) = 2\pi b(\sx_\one-\sx_\two) 
 + i b \partial_t w_b(t)\Bigm|_{t=\spp_2} \, \frac{1}{w_b(\spp_2)}
\end{aligned}
\end{equation}
and, analogously, 
$\log(\sf_b ) = 2\pi b(\sx_\two-\sx_\one) +  b \psi_b(\spp_1)$.
Exponentiating these relations, we obtain~(\ref{expEF}).
\end{proof}

\begin{rem}
Alternatively, eqs.~(\ref{expEF}) can be derived from
(\ref{EFK2}) with the help of the Baker-Campbell-Hausdorff 
formula. Observe also that
$(\log\se_b + \log\sf_b) = b(\psi_b(\spp_1)+\psi_b(\spp_2))$
commutes with $\SC$ and~$\SH$. 
\end{rem}

We may now give the proof of Lemma \ref{codual}.
In Lemma~\ref{uvg}, let us choose 
$A= (2\pi(\sx_\one-\sx_\two) +  \psi_b(\spp_2))\ot 1
 + 1 \ot \frac{\pi}{2} (\spp_2 -\spp_1)$ and
$B= 1\ot (2\pi(\sx_\one-\sx_\two) +  \psi_b(\spp_2))
 + \frac{\pi}{2} (\spp_1 -\spp_2) \ot 1$.
In view of (\ref{expEF}) this implies the identification
$u=\se_b \ot \SK_b$ and $v=\SK_b^{-1} \ot \se_b$.
Therefore, we have
\begin{equation}\nonumber
\begin{aligned}
 \Delta( \se_\frac{1}{b}) & \={efdual} 
 \Delta \bigl( (\se_b)^\frac{1}{b^2} \bigr) =
 \bigl(\Delta (\se_b) \bigr)^\frac{1}{b^2} \={De}
 \bigl(\se_b \ot \SK_b + \SK_b^{-1} \ot \se_b 
   \bigr)^\frac{1}{b^2} \\
 & \={qF} \bigl(\se_b \ot \SK_b \bigr)^\frac{1}{b^2} + 
 \bigl( \SK_b^{-1} \ot \se_b \bigr)^\frac{1}{b^2}
 \={efdual} \se_\frac{1}{b} \ot \SK_\frac{1}{b} +
  \SK_\frac{1}{b}^{-1} \ot \se_\frac{1}{b} \,.
\end{aligned}
\end{equation}
The relation for $\tilde{\SF}$ in (\ref{ded}) is proven similarly.
\qed

\begin{lem}\label{eft}
Let the powers $E^{\alpha}$ and $F^{\alpha}$ with 
$\Im \alpha \neq 0$ be defined
on $\CH$ in the sense of (\ref{EFphi}). Then we have
\begin{equation}\label{EFt}
 [ \SE , \SF^{\alpha} ] = 
  [\alpha]_q \, [ 2 \SH + \alpha -1 ]_q \, \SF^{\alpha-1} \,,
 \quad 
 [ \SE^{\alpha} , \SF ] = 
  [\alpha]_q \, [ 2 \SH - \alpha +1 ]_q \, \SE^{\alpha-1} \,,
\end{equation}
where the $q$-numbers are defined as in Theorem~\ref{T2}.
\end{lem}
\begin{proof}
\begin{equation*}
\begin{aligned}
 {} [ \SE , \SF^{\alpha} ] & \={EFw}  
 \Bigl[ w_b(\spp_2) \, 
 \frac{e^{2\pi b(\sx_\one-\sx_\two)}}{2\sin\pi b^2} \, 
 w_b(-\spp_2) \, , \,  w_b(\spp_1) \, 
 \frac{e^{2\pi \alpha b(\sx_\two-\sx_\one)}}{(2\sin\pi b^2)^\alpha} 
 \, w_b(-\spp_1) \Bigr] \\
 & \={fe}  \Bigl( \frac{w_b(\spp_1 + ib)\, w_b(\spp_2)}%
 {w_b(\spp_1)\, w_b(\spp_2 - i b)} -
  \frac{w_b(\spp_1 + (1-\alpha) ib)\, w_b(\spp_2 + \alpha ib)}%
 {w_b(\spp_1 - \alpha ib)\, w_b(\spp_2 + (\alpha-1) ib)} \Bigr) \\
 & \phantom{===} \times \, w_b(\spp_1) \, 
 \frac{e^{2\pi (\alpha-1) b(\sx_\two-\sx_\one)}}%
 {(2\sin\pi b^2)^{\alpha+1}}  \, w_b(-\spp_1) \\
 & \={wfunrel} 
 \Bigl( \bigl[\fr{Q}{2b} - \fr{i \spp_1}{b}\bigr]_q \,
 \bigl[\fr{Q}{2b} + \fr{i \spp_2}{b}\bigr]_q -
 \bigl[\fr{Q}{2b} -\alpha - \fr{i \spp_1}{b}\bigr]_q \,
 \bigl[\fr{Q}{2b} -\alpha + \fr{i \spp_2}{b}\bigr]_q 
 \Bigr) \, \SF^{\alpha-1} \\
 &= [\alpha]_q \, \bigl[\fr{Q}{b} -\alpha + 
   \fr{i}{b}(\spp_2 - \spp_1) \bigr]_q \, \SF^{\alpha-1}
 = [\alpha]_q \, \bigl[ \fr{1}{ib}(\spp_2 - \spp_1) 
   +\alpha -1 \bigr]_q \, \SF^{\alpha-1} \\
 & \={EFK2}  [\alpha]_q \, \bigl[ 2 \SH
   +\alpha -1 \bigr]_q \, \SF^{\alpha-1} \,.
\end{aligned}
\end{equation*}
In the fourth line we used the definition of 
$q$-number; in the fifth line we used the identities 
$[x]_q [y]_q - [x-\alpha]_q [y-\alpha]_q = [\alpha]_q 
 [x+y-\alpha]_q$ and $[t+b^{-2}]_q = -[t]_q$.
The second formula in (\ref{EFt}) is derived analogously.
\end{proof}

\begin{lem}\label{Qexp}
Let $A$ and $B$ be self-adjoint operators such that
$[A,B] = 2 \pi i$. Then for the function $g_b(x)$
defined in (\ref{gG}) we have
\begin{eqnarray}
 \label{qexp}
 g_b(u) \ g_b(v) &=& g_b(u+v) \,, \\
 \label{pent}
 g_b(v) \ g_b(u) &=& g_b(u) \ g_b(q^{-1} uv) \ g_b(v) \,,
\end{eqnarray}
where $u=e^{bA}$, $v=e^{bB}$ and $q=e^{i \pi b^2}$.
Furthermore, (\ref{qexp}) $\Leftrightarrow$ (\ref{pent}).
\end{lem}
In the literature, eqs.~(\ref{qexp}) and (\ref{pent}) 
are often referred to 
as the quantum exponential and the quantum pentagon
relations. They also hold for the function
$s_q(x) = \prod_{n=0}^\infty (1+xq^{2n+1})$ which is
the compact counterpart of $g_b(x)$. For $s_q(x)$, 
the quantum exponential relation has been known 
since long time \cite{Sch} and the quantum pentagon
relation was found in~\cite{FV}.

Since (\ref{qexp}) and (\ref{pent}) are equivalent,
it suffices to prove one of them. Proofs of 
the quantum pentagon relation were given in 
\cite{FVK} and~\cite{W2}. Nevertheless, we find it
instructive to give another proof of the quantum 
exponential relation in Appendix B since it will allow 
us to introduce the notion of $b$-binomial coefficients.

\subsection{Alternative representations of the R-operator}

\begin{lem}\label{Rint}
$\SR$ and ${\SR}^{-1}$ may be decomposed into 
powers of\, $\SE\ot\SF$ as follows:
\begin{eqnarray}
 \SR &=& b \, \int\limits_{\BR}dt\;
\frac{e^{-\pi i b^2 t^2}}{G_b(Q+ibt)}\;
  q^{\SH\ot\SH} \bigl( 4 (\sin \pi b^2)^2 \, 
  \SE\ot\SF\bigr)^{it}\, q^{\SH\ot\SH} \,, 
\label{Rint1} \\
 {\SR}^{-1} &=& b \, \int\limits_{\BR}dt\;
\frac{e^{-\pi b Q t}}{G_b(Q+ibt)}\;
 q^{-\SH\ot\SH} \bigl( 4 (\sin \pi b^2)^2 \, 
 \SE\ot\SF\bigr)^{it}\, q^{-\SH\ot\SH} \,, 
 \label{Rint2}
\end{eqnarray}
where the integration contour goes above the pole at $t=0$.
\end{lem}
\begin{proof}
By Lemma 15 in \cite{PT} (see also 
\cite{FVK}\cite{K}) we have: 
\begin{equation}\label{bbeta}
\int\nolimits_\BR d\tau \;e^{- 2\pi \tau \beta} \,
\frac{G_b(\alpha + i\tau)}{G_b(Q+ i \tau)} \,=\,
\frac{G_b(\alpha) \, G_b(\beta)}{G_b(\alpha+\beta)} \,.
\end{equation}
The function $\frac{1}{G_b(Q+i\tau)}$ has a pole at $\tau=0$
and is analytic in the upper half-plane. The integration 
contour in (\ref{bbeta}) goes above this pole.

Considering the asymptotics of (\ref{bbeta}) for
$\Im\al\ra -\infty$ and $\Im\al\ra +\infty$, 
using the properties (\ref{Grefl}) and (\ref{Gasym}), 
and making a change of variables, we obtain the
following Fourier transformation formulae 
\begin{eqnarray} 
\label{GG1}
&& b \! \int\limits_{\BR + i 0} \! dt \, e^{2\pi i btr}\,
 \frac{e^{-\pi i b^2t^2}}{G_b(Q+ibt)}  = 
 \frac{\overline{\zeta}_b}{G_b\bigl(\fr{Q}{2}-ir\bigr)} 
 = g_b(e^{2\pi b r}) \,, \\
\label{GG2}
&& b \! \int\limits_{\BR + i 0} \! dt \, e^{2\pi i b t r }\,
 \frac{ e^{-\pi b Q t} }{G_b(Q+ibt)} = 
 \zeta_b \, G_b\bigl(\fr{Q}{2}-ir\bigr) =
 \bigl( g_b(e^{2\pi b r}) \bigr)^{-1} \,.
\end{eqnarray}
Lemma  \ref{Rint} follows if we put here
$r=\frac{1}{2\pi b} \log (4 (\sin \pi b^2)^2 \SE \ot \SF)$
and compare the result with the definition (\ref{Rdef}).
\end{proof}

We now come to the proof of Theorem~\ref{Tint}.
Consider the product representation of the $\SR$-operator
of Lemma~\ref{Rprod} projected to 
$\CP_{s_{\two}}\ot \CP_{s_{\one}}$
by means of the reduction described
in Subsection~\ref{CH},
\begin{eqnarray} \label{Rss}
 \SR_{s_\two s_\one} & \={Rpr} & 
 q^{\SH_\two \ot \SH_\one} \, 
 \bigl( w_b(s_\two - \spp)  \ot w_b(s_\one + \spp) \bigr) \; 
 g_b \bigl( e^{2\pi b \sx}  \ot e^{-2\pi b \sx} \bigr) \\
 \nonumber
 && \hspace*{5mm} \times
 \bigl( w_b(s_\two - \spp) \ot w_b(s_\one + \spp) \bigr)^{-1}\,
 q^{\SH_\two \ot \SH_\one } \\
 \nonumber
 & \={GG1} & q^{\SH_\two \ot \SH_\one} \, 
  \bigl( w_b(s_\two - \spp)  \ot w_b(s_\one + \spp) \bigr) \,
 \int\nolimits_{\BR}d\tau \,  
 \frac{e^{-\pi i \tau^2} e^{2\pi i \tau \sx} \otimes 
 e^{-2\pi i \tau \sx}}{G_b(Q+i\tau)}\, \\
 \nonumber
 && \hspace*{5mm} \times \bigl( w_b(s_\two - \spp) \otimes 
 w_b(s_\one + \spp) \bigr)^{-1}\,
 q^{\SH_\two \ot \SH_\one }
\end{eqnarray}
It is now easy to compute the ``matrix elements''  of 
$\SR_{s_\two s_\one}$ on the states 
$| k_2, k_1 \rangle = | k_2 \rangle \otimes | k_1 \rangle$,
where $| k \rangle \equiv e^{2\pi i x k}$.
Taking into account that $\spp | k \rangle = k | k \rangle$
and $\langle k | k' \rangle = \delta(k' - k)$, we find
\begin{eqnarray}
 \label{Skk1}
 &&  \tilde{\SR}_{s_2s_1}  (k_2,k_1| k'_2,k'_1) \,=\, 
 \langle k_2, k_1 | \SR_{s_2s_1} | k'_2, k'_1 \rangle \\
\nonumber
 && = \int\limits_{\BR}d\tau \,
 \frac{e^{-\pi i (\tau^2 + k_1 k_2 + k'_1 k'_2)}}{G_b(Q+i\tau)}\,
 \frac{w_b(s_1 + k_1)}{w_b(s_1 + k'_1)} \,
 \frac{w_b(s_2 - k_2)}{w_b(s_2 - k'_2)} \,
 \delta(k'_2 - k_2 + \tau) \,
 \delta(k'_1 - k_1 - \tau) 
\end{eqnarray} 
which gives us the kernel (\ref{ker2})
of the ``momentum'' representation in Theorem~\ref{Tint}.

The kernel of the ``coordinate'' representation (\ref{io1})
can be obtained as a Fourier transform of (\ref{Skk1}):
\begin{equation} \nonumber
\begin{aligned}
 \SR_{s_2s_1} & (x_2,x_1|x'_2,x'_1) =
 \int_\BR dk'_2 dk'_1 dk_2 dk_1 \, 
 e^{2\pi i(x_2 k_2 + x_1 k_1 - x'_2 k'_2 - x'_1 k'_1) } \,
 \tilde{\SR}_{s_2s_1}(k_2,k_1|k'_2,k'_1) \\
 &= \int_\BR d\tau dk'_2 dk'_1 \,
 \frac{e^{\pi i \left( \tau (k'_2 - k'_1) - 2 k'_1 k'_2 \right)}}%
 {G_b(Q+i\tau)}\, 
 e^{2 \pi i \left( \tau(x_2 - x_1) + k'_1 (x_1 - x'_1)
 + k'_2 (x_2 - x'_2) \right)} \\
 &\ \hspace*{5mm} \times 
 \frac{w_b(s_1 + k'_1 - \tau)}{w_b(s_1 + k'_1)} \,
 \frac{w_b(s_2 - k'_2 - \tau)}{w_b(s_2 - k'_2)} \,.
\end{aligned}
\end{equation}
The remaining integrations are performed by
using relation (\ref{bbeta}) three times. The
result of this straightforward but tedious 
calculation is given by~(\ref{ker1}).

\section{Proofs of the main results}
\setcounter{equation}{0}

\subsection{Proof of Theorem \ref{coass}: Co-associativity}

First, it is straightforward to write out $(\id\ot\De)\circ\De(\SX)$
and ${(\De\ot\id)}\circ{\De(\SX)}$
in terms of the Clebsch-Gordan maps $\SC[s_3|s_2,s_1]$:
\[\begin{aligned}
(\pi_{s_3}\ot\pi_{s_2}\ot\pi_{s_1})\circ
(\id\ot\De)\circ\De(\SX)\;=\;& \int_{\BR^+}d{\rm m}(s_4)d{\rm m}(s_{21})\;
\SC_{3(21)}^{\dagger}(s_{21})\cdot\SX \cdot\SC_{3(21)}(s_{21})\\
(\pi_{s_3}\ot\pi_{s_2}\ot\pi_{s_1})\circ
(\De\ot\id)}\circ{\De(\SX)\;=\;& \int_{\BR^+}d{\rm m}(s_4)d{\rm m}(s_{32})\;
\SC_{(32)1}^{\dagger}(s_{32})\cdot\SX \cdot\SC_{(32)1}(s_{32}),
\end{aligned}
\]
where we have introduced
\begin{equation}\begin{aligned}
\SC_{3(21)}(s_{21})\;\equiv\;& 
\SC[s_4|s_3,s_{21}]\cdot\big(\id\ot\SC[s_{21}|s_2,s_1]\big),\\
\SC_{(32)1}(s_{32})\;\equiv\;& 
\SC[s_4|s_{32},s_{1}]\cdot\big(\SC[s_{32}|s_3,s_2]\ot\id\big).
\end{aligned}\end{equation}
Proposition 7 in \cite{PT} is equivalent to 
\[
\SC_{3(21)}(s_{21})\;=\;\int_{\BR^+}ds_{32}\;
\big\{\begin{smallmatrix} s_1 & s_2 & s_{21} \\
s_3 & s_4 & s_{32}\end{smallmatrix}\big\}\,\SC_{(32)1}(s_{32}),
\]
where $\big\{\begin{smallmatrix} s_1 & s_2 & s_{21} \\
s_3 & s_4 & s_{32}\end{smallmatrix}\big\}_b$ are the b-Racah-Wigner 
coefficients introduced in \cite{PT}.
It follows that 
\[ 
\begin{aligned}
(  \pi_{s_3}\ot & \pi_{s_2}\ot\pi_{s_1})\circ 
(\id\ot\De)\circ\De(\SX)\;=\;\\
& =\;\int\limits_{\BR^+}
d{\rm m}(s_4)d{\rm m}(s_{21})\int\limits_{\BR^+}d{\rm m}(s_{32})
d{\rm m}(s_{32}') \;
\big\{\begin{smallmatrix} s_1 & s_2 & s_{21} \\
s_3 & s_4 & s_{32}' \end{smallmatrix}\big\}^{*}_b
\big\{\begin{smallmatrix} s_1 & s_2 & s_{21} \\
s_3 & s_4 & s_{32}\end{smallmatrix}\big\}_b^{}
\cdot \\ & \hspace{7cm}\cdot
\SC_{(32)1}^{\dagger}(s_{32})\cdot\SX \cdot\SC_{(32)1}(s_{32}).
\end{aligned}
\]
Exchanging the integrations over $s_{21}$ and $s_{32}$, $s_{32}'$,
and using formula (89) from \cite{PT},
\[ 
\int\limits_{\BR^+}d{\rm m}(s_{21})\;
\big\{\begin{smallmatrix} s_1 & s_2 & s_{21} \\
s_3 & s_4 & s_{32}' \end{smallmatrix}\big\}^{*}_b
\big\{\begin{smallmatrix} s_1 & s_2 & s_{21} \\
s_3 & s_4 & s_{32}\end{smallmatrix}\big\}^{}_b
={\rm m}(s_{32})
\de(s_{32}-s_{32}'),
\] 
yields the claim.

\subsection{Proof of Theorem 2: Invariance of the Haar-measure}

We shall consider the left invariant Haar measure $\sh_{\rm l}$ only,
the proof for the case of $\sh_{\rm r}$ being completely analogous.

A few preparations are in order.
The elements of $\CA_{\sh}^{\rm l}$ can be represented as integral operators:
If a vector $\psi\in\CH$ is realized by a function $\psi(k,s)$, then
\[
(\SO\psi)(k,s)\;=\;\int_\BR dk'\;K_{\SO}(k,k'|s)\psi(k',s).
\]
In terms of the kernel $K_{\SO}(k,k'|s)$ one may write the
defintion of $\sh_{\rm l}$ as
\begin{equation}\label{haarrep}
\sh_{\rm l}(\SO)\;=\;\int_{\BR^+}d{\rm m}(s)\int_{\BR}dk\; e^{-2\pi Qk}
K_{\SO}(k,k|s).
\end{equation}
The distributional matrix elements of an operator $\SO\in\CA$ are
always of the form
\begin{equation}
\bra s,k| \SO|s',k'\ket\,\equiv\,
\de(s-s')
\bra\!\!\bra k|\SO|k'\ket\!\!\ket_{s}^{}.
\end{equation}
By using an analogous notation for operators in $\CA\ot\CA$ 
one may represent the distributional matrix elements of $\De(\SO)$ as
\begin{equation}\label{corep}\begin{aligned}
\bra\!\!\bra k_2,k_1| & \De(\SO)|k_2',k_1'\ket\!\!\ket_{s_2s_1}\\
& =\int_{\BR^+}d{\rm m(s_3)}\int_{\BR}dk_3dk_3'\;K_{\SO}(k_3,k_3'|s_3)
\Big(\CGC{s_3}{k_3}{s_2}{k_2}{s_1}{k_1}\Big)^*
\CGC{s_3}{k_3'}{s_2}{k_2'}{s_1}{k_1'} .
\end{aligned}\end{equation}

In order to make the justification for the following 
manipulations more transparent, we are going
to employ the following regularization for the distributions involved:
\[ 
\CGC{s_3}{k_3}{s_2}{k_2}{s_1}{k_1}\;=\;\lim_{\ep\downarrow 0}\;
\CGC{s_3}{k_3}{s_2}{k_2}{s_1}{k_1}_{\ep},\qquad
\CGC{s_3}{k_3}{s_2}{k_2}{s_1}{k_1}_{\ep}\;=\;
e^{-\ep\sum_{i=1}^3|k_i|}\de_{\ep}(k_3-k_2-k_1)
\CGCtr{s_3}{s_2}{k_2}{s_1}{k_1},
\]
where $\de_{\ep}(x)=\de_{\ep}(-x)$ is a 
symmetric regularization of the delta-distribution.
Let us furthermore note that it suffices to check the invariance
property on a dense subset $\CT$ of the domain of $\sh_{\rm l}$. 
Consider the matrix element
\[
\bra\!\!\bra \psi_2,k_1|\De(\SO)|\psi_2',k_1'\ket\!\!\ket_{s_2s_1}^{}
\; :=\; 
\int_{\BR}dk_2^{}dk_2'\;\psi(k_2)\psi'(k_2')
\bra\!\!\bra k_2,k_1|\De(\SO)|k_2',k_1'\ket\!\!\ket_{s_2s_1}^{},
\]
where $\psi_2$, $\psi_2'$ are smooth functions with compact support.
Assuming 
that $K_{\SO}(k,k'|s)$
has exponential decay w.r.t. $k$ and $k'$
it is not difficult to show that the matrix element
$\bra\!\!\bra \psi_2,k_1|\De(\SO)|\psi_2',k_1'\ket\!\!\ket_{s_2s_1}^{}$
will also have exponential decay w.r.t. $k_1$, $k_1'$ that can be
made as large as one likes by choosing the subset $\CT\subset\CA_h$
appropriately.

Combining \rf{haarrep} and \rf{corep} 
leads to the following representation for 
the distributional matrix elements of $(\id\ot\sh_{\rm l})\De(\SO)$.
\begin{equation}\label{inv_m_el}
\begin{aligned}
\bra\!\!\bra k_2|(\id\ot\sh_{\rm l})\De(\SO)|k_2'\ket\!\!\ket_{s_2}^{}=
& \int_{\BR^+} d{\rm m(s_3)}  \int_{\BR}dk_3dk_3'\;K_{\SO}(k_3,k_3'|s_3)\\
&\ti \int_{\BR^+}d{\rm m(s_1)}\int_{\BR}dk_1\; e^{-2\pi Qk_1}
\Big(\CGC{s_3}{k_3}{s_2}{k_2}{s_1}{k_1}\Big)^*
\CGC{s_3}{k_3'}{s_2}{k_2'}{s_1}{k_1} .
\end{aligned}\end{equation}
We are going to use the following result: 
\begin{propn}
The following equation holds as an identity
between tempered distributions:
\begin{equation}\label{dist_id}
\begin{aligned} 
e^{2\pi k_3Q}
\int_{\BR^+}d{\rm m(s_1)}\int_{\BR}dk_1\; e^{-2\pi Qk_1}
\Big(\CGC{s_3}{k_3}{s_2}{k_2}{s_1}{k_1}\Big)^* & 
\CGC{s_3}{k_3'}{s_2}{k_2'}{s_1}{k_1}
\;=\;\\
& =\; \de(k_3'-k_3)\de(k_2'-k_2).
\end{aligned} 
\end{equation}
\end{propn}
\begin{proof} 
The proof of the proposition will be based on the 
following important symmetries of the Clebsch-Gordan kernel:
\begin{lem}
The Clebsch-Gordan kernel $\CGC{s_3}{k_3}{s_2}{k_2}{s_1}{k_1}_{\ep}$
has the following symmetries.
\begin{equation}\label{CGCksymm}\begin{aligned}
\Big(\CGC{s_3}{k_3}{s_2}{k_2}{s_1}{k_1}_{\ep}\Big)^*\;=\;
e^{+\pi Q( k_1-k_3)}e^{-\pi i h_{s_2}}
\CGC{s_1}{-k_1}{-s_2}{-k_2}{s_3}{-k_3}_{\ep},\\
\Big(\CGC{s_3}{k_3}{s_2}{k_2}{s_1}{k_1}_{\ep}\Big)^*\;=\;
e^{-\pi Q( k_2-k_3)}e^{-\pi i h_{s_1}}
\CGC{s_2}{-k_2}{s_3}{-k_3}{-s_1}{-k_1}_{\ep}.
\end{aligned}\end{equation} 
\end{lem}
\begin{proof}
One may verify directly that for $s_i\in\BR$, $i=1,2,3$, 
\begin{equation}\begin{aligned}\label{CGCxsymm}
\Big(\CGC{s_3}{x_3}{s_2}{x_2}{s_1}{x_1}\Big)^*\;=\;
e^{-\pi i h_{s_2}}\CGCb{s_1}{x_1^*-i\frac{Q}{2}}
{-s_2}{x_2^*}{s_3}{x_3^*-i\frac{Q}{2}},\\
\Big(\CGC{s_3}{x_3}{s_2}{x_2}{s_1}{x_1}\Big)^*\;=\;
e^{-\pi i h_{s_1}}\CGCb{s_2}{x_2^*+i\frac{Q}{2}}
{s_3}{x_3^*+i\frac{Q}{2}}{-s_1}{x_1^*}.
\end{aligned}\end{equation}
The Lemma follows by taking the Fourier-transformation of \rf{CGCxsymm},
taking into account that our regularization is compatible with the symmetry
\rf{CGCksymm}.
\end{proof}
With the help of equation \rf{CGCksymm} we may rewrite the 
left hand side of \rf{dist_id} as follows:
\begin{equation}
\lim_{\ep\downarrow 0}\; e^{\pi (k_3-k'_3)Q}
\int_{\BR^+}d{\rm m(s_1)}\int_{\BR}dk_1\; 
\Big(\CGC{s_1}{-k_1}{-s_2}{-k_2}{s_3}{-k_3}_{\ep}\Big)^*
\CGC{s_1}{-k_1}{-s_2}{-k_2'}{s_3}{-k_3'}_{\ep}  
\end{equation}
The proposition now follows by using the Fourier-transform of
\cite[Corollary 1]{PT}.
\end{proof}
Inserting \rf{dist_id} into \rf{inv_m_el} yields
\[ 
 \bra\!\!\bra k_2|(\id\ot\sh_{\rm l})
 \De(\SO)|k_2'\ket\!\!\ket_{s_2}^{}\;=\;
\de(k_2-k_2')\int_{\BR^+}d{\rm m}(s_3)\int_{\BR}dk_3\; e^{-2\pi Qk_3}
K_{\SO}(k_3,k_3|s_3).
\]
Recognizing the definition of the Haar-measure on the right hand 
side completes the proof of the left invariance property
of $\sh_{\rm l}$.

To prove the property (ii) in Theorem~2 we observe that
definition (\ref{Hm}) can be rewritten as 
\[
 \sh_{\rm l}(\SO) \;=\;
 \int_{0}^{\infty}d{\rm m}(s)\;{\rm Tr}(\SK^2 \tilde{\SK}^2 \SO_s) \,.
\]
Now it is straightforward to verify the first formula
in (\ref{Hinv2}) for $\SX=\SE, \SF, \SK$ using the definition of
the adjoint action, the relations (\ref{def}), and the cyclicity
of trace. For instance, 
${\rm Tr}(\SK^2 \tilde{\SK}^2 {\rm ad^l_\SK}(\SO_s))=
 {\rm Tr}(\SK^2 \tilde{\SK}^2 \SK \SO_s \SK^{-1}) =
 {\rm Tr}(\SK^2 \tilde{\SK}^2 \SO_s )$,
${\rm Tr}(\SK^2 \tilde{\SK}^2 {\rm ad^l_\SE}(\SO_s))=
 {\rm Tr}(\SK^2 \tilde{\SK}^2 (\SE \SO_s \SK^{-1} -
 q \SK^{-1} \SO_s \SE ))=
 {\rm Tr}((\SK \SE  - 
 q \SE \SK ) \tilde{\SK}^2 \SO_s ) = 0$.
Further, we notice that $\sh_{\rm l}( {\rm ad^l_{\SX\SY}}(\SO))=
 \sh_{\rm l}( {\rm ad^l_\SX} ({\rm ad^l_\SY}(\SO)) )=
 \sh_{\rm l}(\SO)\, \epsilon(\SX) \, \epsilon(\SY) =
 \sh_{\rm l}(\SO)\, \epsilon(\SX \SY)$. Together with
the linearity of trace this implies that (\ref{Hinv2})
extends to any element of $U_q(\fsl(2,\BR))$.

\subsection{Proof of Theorem~\ref{T1}} 

Let us adopt the following notations: 
$\SX_1 \equiv \SX \otimes 1 \otimes 1$,
$\SX_2 \equiv 1 \otimes \SX \otimes 1$, and 
$\SX_3 \equiv 1 \otimes 1 \otimes \SX$.
\newline
{\it\bf Property (i)}\ \  First, we compute
with the help of Lemma \ref{eft}:
\begin{eqnarray}
 &&  \label{ia}
 q^{\SH\ot\SH} \, (\SE_1 \SF_2)^{i t} \, q^{\SH\ot\SH} \,
 \bigl( \SK_1^{-1} \SE_2 \bigr) - 
 \bigl( \SK_1 \SE_2 \bigr) \, 
 q^{\SH\ot\SH} \, (\SE_1 \SF_2)^{i t} \, q^{\SH\ot\SH} \\
 && \nonumber
 = q^{\SH\ot\SH} \, \bigl[ (\SE_1 \SF_2)^{i t} , 
  \SE_2 \bigr]  \, q^{\SH\ot\SH} 
 \={EFt} - q^{\SH\ot\SH} \Bigl( [i t]_q [2 \SH_2 + it -1]_q \,
 \SE_1^{i t} \, \SF_2^{i t-1}  \Bigr) q^{\SH\ot\SH} \,.
\end{eqnarray}
Next, we find
\begin{eqnarray}
 && \label{ib}
  q^{\SH\ot\SH} \, (\SE_1 \SF_2)^{i t} \, q^{\SH\ot\SH} \,
 \bigl( \SE_1 \SK_2 \bigr) - 
 \bigl( \SE_1 \SK_2^{-1} \bigr) \, 
 q^{\SH\ot\SH} \, (\SE_1 \SF_2)^{i t} \, q^{\SH\ot\SH} \\
  && \nonumber 
 = q^{\SH\ot\SH} \, \Bigl( (\SE_1 \SF_2)^{i t} \, 
 \SE_1 \SK_2^2 - \SE_1 \SK_2^{-2} (\SE_1 \SF_2)^{i t}
 \Bigr) \, q^{\SH\ot\SH} \\
 && \nonumber  
 =  q^{\SH\ot\SH} \Bigl( 
 (2 i\, q^{i t} \, \sin \pi b^2) \, [2 \SH_2 + i t]_q  
  \, \SE_1^{i t+1} \, \SF_2^{i t} \Bigr) \, q^{\SH\ot\SH} \,.
\end{eqnarray}
Let us write down the integral representation (\ref{Rint1})
of $\SR$ in the following form
\begin{equation}\label{RI}
 \SR \,=\, \int\nolimits_{\BR}dt\, \rho(t) \, 
 q^{\SH\ot\SH} \bigl(\SE\ot\SF\bigr)^{it}\, q^{\SH\ot\SH} \,, 
\end{equation}
where $\rho(t) \equiv b 
 \frac{e^{-\pi i b^2 t^2}}{G_b(Q+ibt)} (2\sin\pi b^2)^{2it}$.
Observe that (\ref{Gfunrel}) implies that $\rho(t)$ 
satisfies the following functional equation:
\begin{equation}\label{req}
 [it +1]_q \, \rho(t-i) = 
 ( 2 i\, q^{i t} \, \sin \pi b^2 ) \, \rho(t) \,.
\end{equation}
Adding (\ref{ia}) with (\ref{ib}), we derive
\begin{equation*}
\begin{aligned}
 \SR \, \Delta(\SE) &- \Delta'(\SE) \, \SR =
  \SR \, (  \SE_1 \SK_2 + \SK_1^{-1} \SE_2 ) -
  ( \SE_1 \SK_2^{-1}  + \SK_1 \SE_2 ) \, \SR \\
 &  \={RI} \int_{\BR}dt\, q^{\SH\ot\SH} \, \rho(t) \, \Bigl(
 (2 i \, q^{i t} \, \sin \pi b^2) \, [2 \SH_2 + i t]_q  \,
 \SE_1^{i t+1} \, \SF_2^{i t} \\
 & \qquad\qquad\qquad\qquad
 - [i t]_q [2 \SH_2 + it -1]_q \,
 \SE_1^{i t} \, \SF_2^{i t-1} \Bigr) \, q^{\SH\ot\SH} \\
 & = \int_{\BR}dt\, q^{\SH\ot\SH} \, \Bigl(
 ( 2 i \, q^{i t} \, \sin \pi b^2 ) \, \rho(t)  \\
 & \qquad\qquad\qquad\qquad
 - [i t+1]_q \, \rho(t -i)  \Bigr) [2 \SH_2 + it]_q \, 
 \SE_1^{i t+1} \, \SF_2^{i t}  \, q^{\SH\ot\SH} \={req} 0 \,.
\end{aligned}
\end{equation*}
Thus, we have proven (\ref{i}) for $\SX=\SE$. 
The proof for $\SF$ goes along the same lines
with the help of the second formula in~(\ref{EFt}).
And for $\SK$ the proof is trivial because 
$\Delta(\SK)$ commutes with $(\SE \ot \SF)^{i t}$.
\\[1mm]
{\it\bf Property (ii)}\ \
Recall that the rescaled generators $\se_b$ and 
$\sf_b$ were introduced before Lemma~\ref{Rprod}.
To prove the first formula in (\ref{ii}), we use 
the quantum exponential relation (\ref{qexp}) {}from 
Lemma~\ref{Qexp} with identification
\hbox{$u= \se_1 \SK^{-1}_2 \sf_3$} 
and \hbox{$v= \se_1 \sf_2 \SK_3$},
\begin{eqnarray*} 
 (\id\ot\De)\SR & \={Rdef} &  
 (\id\ot\De) \bigl( q^{ \SH_1 \SH_2} \,
 g_b ( \se_1 \sf_2) \, q^{ \SH_1 \SH_2} \bigr) \\
 & \={De} & q^{ \SH_1 \SH_2 + \SH_1 \SH_3} \,
 g_b \bigl( \se_1 \sf_2 \SK_3 + 
  \se_1 \SK^{-1}_2 \sf_3  \bigr)
 \, q^{ \SH_1 \SH_2 + \SH_1 \SH_3 }   \\ 
 & \={qexp} & q^{ \SH_1 \SH_2 + \SH_1 \SH_3 } \, 
 g_b( \se_1 \SK^{-1}_2 \sf_3 ) \,
 g_b( \se_1 \sf_2 \SK_3 ) \,
 q^{ \SH_1 \SH_2 + \SH_1 \SH_3 } \\
 & \={def} & q^{ \SH_1 \SH_3} \, g_b( \se_1 \sf_3 ) \,
 q^{ \SH_1 \SH_3} \, q^{ \SH_1 \SH_2} \, g_b( \se_1 \sf_2 ) 
 \, q^{ \SH_1 \SH_2} =  \SR_{13}\SR_{12} \,. 
\end{eqnarray*}
The second formula in (\ref{ii}) is proved in the same way. 
\\[1mm]
\noindent
{\it\bf Property (iii)}\ \
First, we derive 
\begin{eqnarray}
 && \nonumber
(\sigma \ot \id) \,  \bigl( 
 q^{ \SH_1 \SH_2} (\SE_1 \SF_2 )^{i t}
 q^{ \SH_1 \SH_2} \bigr)  
 = (\sigma \ot \id) \, \bigl(
 q^{ \SH_1 \SH_2} (\SE_1)^{i t}
 q^{ \SH_1 (\SH_2 + i t)} \bigr) (\SF_2 )^{i t} \\
 && \nonumber 
 = q^{- \SH_1 (\SH_2+ i t)} 
  (-q \SE_1)^{i t}  q^{- \SH_1 \SH_2} (\SF_2)^{i t} 
  = q^{- \SH_1 \SH_2} 
 q^{- i t \SH_1} (-q \SE_1)^{i t}  
 q^{ i t \SH_1} (\SF_2)^{i t} q^{- \SH_1 \SH_2} \\
&&  = e^{ \pi b (i b t^2 - Q t )}
 q^{- \SH_1 \SH_2} (\SE_1 \SF_2 )^{i t}
 q^{- \SH_1 \SH_2} \,. \label{sR}
\end{eqnarray}
This means that, acting with $(\sigma \ot \id)$
on the r.h.s.\ of (\ref{Rint1}), we obtain the r.h.s.\ 
of (\ref{Rint2}). Thus, we have proven the first formula 
in (\ref{ii}). The second formula is verified analogously. 
Finally, acting with $(\id \ot \sigma)$ on the last line in
(\ref{sR}) and performing similar manipulations,
we find that  
\begin{equation*}
 (\id \ot \sigma) \bigl( e^{ \pi b (i b t^2 - Q t )}
 q^{- \SH_1 \SH_2} (\SE_1 \SF_2 )^{i t}
 q^{- \SH_1 \SH_2} \bigr) \,=\,
 q^{ \SH_1 \SH_2} (\SE_1 \SF_2 )^{i t} q^{ \SH_1 \SH_2} 
\end{equation*}
which together with (\ref{sR}) implies the last
formula in~(\ref{iii}).

\subsection{Proof of Theorem 5: R-operator the highest weight
representations}

We first need to discuss the analytic continuation of 
\begin{equation}\label{r-act}
\big\bra\, \de_{k_2}\ot\de_{k_1}\, , \SR_{s_2s_1}^{}f\,\big\ket\;=\;
\int_{\BR}dk_2'dk_1'\;\tilde{\SR}(k_2,k_1|k_2',k_2')\,f(k_2',k_2')
\end{equation} 
to the values $k_1=-s_1+i\big(\frac{Q}{2}+nb\big)$, 
$n\in\BZ^{\geq 0}$. 
To begin with, one may trivially perform e.g. the integral over $k_2'$ to
get an expression of the form
\begin{equation}\label{r-act'}
\big\bra\, \de_{k_2}\ot\de_{k_1}\, , \SR_{s_2s_1}^{}f\,\big\ket\;=\;
\int_{\BR}dk_1'\;\tilde{\SR}_{k}(k_1|k_1')\,f(k-k_1',k_1'),
\end{equation}
where $k=k_2+k_1$.
The analytic  continuation of \rf{r-act'} to 
$k_1=-s_1+i\big(\frac{Q}{2}+nb\big)$ can be defined by deforming the
contour of integration over $k_1'$, $\BR$, in \rf{r-act'} into
$\BR+i\big(\frac{Q}{2}+nb)+i0$ plus a sum of small circles around the 
poles from the factor $w_b^{-1}(s_1+k_1')$ in $\tilde{\SR}_{k}(k_1|k_1')$ 
that lie between $\BR$ and 
$ \BR+i\big(\frac{Q}{2}+nb)+i0$.

We are now in the position to take the limit $k_1\ra 
-s_1+i\big(\frac{Q}{2}+nb\big)$. 
The factor $w_b(s_1+k_1)$
that appears in $\tilde{\SR}_{k}(k_1|k_1')$, cf.
\rf{ker2}, 
makes most of the terms vanish except for the
terms from the poles at $k_1'=-s_1+i(\frac{Q}{2}+bn')$, $0\leq n'\leq n$. 
The resulting expression is of the following form.
\begin{equation}
\SR_{s_2s_1}^{t}\;\de_{k_2}\ot  \de_{k_1}=
\sum_{l=0}^{n}G_{l}\,
e^{-\pi i ((k_1-ibl)k_2 + k_1 k'_2)}
\frac{w_{n-l}}{w_n} \,
\frac{w_b(s_2 - k_2)}{w_b(s_2 - k_2')}\bigg|_{k_2'=k_2+ibl}
\end{equation}
where $G_n:=\Res_{x=nb}G_b^{-1}(Q+x)$ and 
$w_n=\Res_{x=i\frac{Q}{2}+inb}w_b^{-1}(x)$.

It remains to calculate the relevant residues. 
It is easy to derive from (\ref{Gdual})-(\ref{Grefl})
that
\begin{equation}\label{Gxx}
  G_b(x) \, G_b(-x) = - \frac{e^{\pi i x^2}}%
 {4 \sin\pi b x \,  \sin\pi b^{-1} x } \ .
\end{equation}
Hence $\lim_{x \rightarrow 0} (x G_b(x))^2 = (2\pi)^{-2}$.
In fact, using the modular property of the Dedekind 
$\eta$-function, it is straightforward to compute the
limit directly for the product representation (\ref{Gprod}) 
(as was done in \cite{Sh}); which yields  
\begin{equation}\label{Glim}
 \lim_{x \rightarrow 0} x \, G_b(x) = \frac{1}{2\pi} \,.
\end{equation}
Hence, taking into account the properties 
(\ref{Gdual})-(\ref{Grefl}), we find that
\begin{equation}\label{Gres}
 \text{ Res}\, \frac{1}{G_b(Q+z)} 
 = -  \frac{1}{2\pi} \, 
 \prod_{k=1}^n (1-q^{2k})^{-1} \,
 \prod_{l=1}^m (1-\tilde{q}^{-2l})^{-1} 
  \quad \text{ at}\ z = nb+mb^{-1} \,,  
\end{equation}
where $n,m \in \BZ^{\geq 0}$ and $\tilde{q}=e^{-i\pi b^{-2}}$.

To complete the proof of Theorem \ref{T2} is now the matter of a 
straightforward calculation using the functional 
relation \rf{wfunrel}, formula \rf{Gres}, as well as \rf{deltatrf}.

\subsection{Proof of Theorem~\ref{T3}}

Let us first note that the left hand side of \rf{Rdiag} 
satisfies the intertwining property
\begin{equation} \label{intertw2}
\SC[s_3|s_1,s_2]\;\SB_{s_{\two}s_{\one}}^{}\circ
(\pi_{s_{\two}}\ot \pi_{s_{\one}})\circ\De(X)\,=\,
\pi_{s_\three}(X)\circ\SC[s_3|s_1,s_2]\;
 \SB_{s_{\two}s_{\one}}^{} \,.
\end{equation}
A unitary operator that maps 
$\CP_{s_{\two}}\ot \CP_{s_{\one}}\ra \CP_{s_{\three}}$ 
and satisfies \rf{intertw} must be proportional to 
$\SC[s_3|s_2,s_1]$. This is a consequence of the analysis 
used to prove Theorem 2 in~\cite{PT}. It follows that there 
{\it exists} a function $\Om(s_3|s_2,s_1)$ such that the statement of 
Theorem~\ref{T3} holds. We are left with the task to calculate
$\Om(s_3|s_2,s_1)$ explicitly.

To this aim let us note that Theorem~\ref{T3} is equivalent
to an identity between meromorphic functions. To write 
this identity down, let us assume that $\CP_{s_{\two}}\ot \CP_{s_{\one}}$
is realized by functions $\tilde{\psi}(k_2,k_1)$. 
$\SC[s_3|s_1,s_2]$ is then realized as an integral operator:
\[
\big(\SC[s_3|s_1,s_2]\tilde{\psi}\big)(k_3)\;=\;\int_{\BR}dk_2dk_1
\;\CGC{s_3}{k_3}{s_2}{k_2}{s_1}{k_1}\tilde{\psi}(k_2,k_1).
\]
The explicit expression for the distributional kernel 
$\CGC{s_3}{k_3}{s_2}{k_2}{s_1}{k_1}$ can be found in Appendix A.
For the moment it will be enough to note that
it can be factorized as
\begin{equation}\label{CGCfact}
\CGC{s_3}{k_3}{s_2}{k_2}{s_1}{k_1}\;=\;
\de(k_1+k_2-k_3)\,\CGCtr{s_3}{s_2}{k_2}{s_1}{k_1},
\end{equation}
where $\CGCtr{s_3}{s_2}{k_2}{s_1}{k_1}$
is a meromorphic function by Lemma 20 of \cite{PT}. 
It is then easy to see that Theorem~\ref{T3}
is equivalent to the identity
\begin{equation}\label{diagcomp}
\int_{\BR}dk_2dk_1\; \CGC{s_3}{k_3}{s_1}{k_1}{s_2}{k_2}\;
\tilde{\SR}_{s_2s_1}(k_2,k_1|k'_2,k'_1)\;=\;
\Om(s_3|s_2,s_1)\CGC{s_3}{k_3}{s_2}{k_2'}{s_1}{k_1'}.
\end{equation}
In order to see that equation \rf{diagcomp} is indeed 
equivalent to an identity between meromorphic functions
let us note that \rf{CGCfact} and \rf{ker2} allow one to 
split off the distributional factors. What remains on the
left hand side is a convolution of two meromorphic functions, so is
itself meromorphic (cf. Lemma 3 in \cite{PT}). 

Let us note that both sides of \rf{diagcomp} have a pole
at $k_2'=-s_2+i\frac{Q}{2}$. In the case of the right hand side this 
is a consequence of Lemma 20 of \cite{PT}. Concerning the left hand side
of \rf{diagcomp} one may as in the proof of Lemma 3 of \cite{PT}
identify the above-mentioned pole as the
consequence of the pinching of the contour of integration by a collision
of two poles of the integrand.
First we have the pole of $\CGC{s_3}{k_3}{s_1}{k_1}{s_2}{k_2}$
at $k_2=-s_2+i\frac{Q}{2}$. Second 
let us note that the factor 
$G_b^{-1}(Q+i(k'_1 - k_1))$ appearing in 
$\tilde{\SR}_{s_2s_1}(k_2,k_1|k'_2,k'_1)$ produces a pole
at $k_2=k_2'$ if one takes into account that 
$\tilde{\SR}_{s_2s_1}(k_2,k_1|k'_2,k'_1)$ has support only for
$k'_1 - k_1=k_2 - k_2'$. The residue of the resulting pole on the
left hand side of \rf{diagcomp} is simply given by the product
of the relevant residues of $\CGC{s_3}{k_3}{s_1}{k_1}{s_2}{k_2}$
and $\tilde{\SR}_{s_2s_1}(k_2,k_1|k'_2,k'_1)$ respectively. 

The equality of the residues of the two sides of
equation \rf{diagcomp} 
implies the following identity:
\begin{equation}
\Res_{k_2=-s_2+i\frac{Q}{2}}
\CGCtr{s_3}{s_1}{k_1}{s_2}{k_2}\; e^{\pi Qk_1}
e^{2\pi is_2k_1}\;=\;
\Om(s_3|s_2,s_1)\Res_{k_2=-s_2+i\frac{Q}{2}}
\CGCtr{s_3}{s_2}{k_2}{s_1}{k_1}
\end{equation}
By evaluating the relevant residues we may therefore 
calculate $\Om(s_3|s_2,s_1)$. 
\begin{lem}
\[ \begin{aligned}
2\pi i\!\!\!\!\!\!\!
\Res_{k_2=-s_2+i\frac{Q}{2}}\CGCtr{s_3}{s_1}{k_1}{s_2}{k_2}
\;=\; & e^{-\frac{\pi i}{2}(h_{s_3}-h_{s_2}-h_{s_1})}e^{-\frac{\pi}{2}Qk_1}
e^{-\pi is_2k_1}
\frac{w_b(k_1-s_1)w_b(s_3+s_1-s_2)}{w_b(k_1+s_3-s_2+i\frac{Q}{2})},\\
2\pi i\!\!\!\!\!\!\!\Res_{k_2=-s_2+i\frac{Q}{2}}\CGCtr{s_3}{s_2}{k_2}{s_1}{k_1}
\;=\; & e^{+\frac{\pi i}{2}(h_{s_3}-h_{s_2}-h_{s_1})}e^{+\frac{\pi}{2}Qk_1}
e^{+\pi is_2k_1}
\frac{w_b(k_1-s_1)w_b(s_3+s_1-s_2)}{w_b(k_1+s_3-s_2+i\frac{Q}{2})}.
\end{aligned} \]
\end{lem}
\begin{proof}
In order to exhibit the singular behavior of 
$\CGCtr{s_3}{s_1}{k_1}{s_2}{k_2}$ near 
 $k_2=-s_2+i\frac{Q}{2}$ one 
may deform the contour of integration in \rf{CGCFT} into the union
of a small circle around the pole of the integrand at $s=0$ and a contour
that separates the pole at $s=0$ from all the other poles
in the {\it upper} half-plane, approaching asymptotically $\pm i\infty$.  
The contribution from the residue of the pole at $s=0$ 
exhibits the pole at
$k_2=-s_2+i\frac{Q}{2}$ explicitly, whereas the rest is nonsingular. 

Similarly, 
to analyze the singular behavior of $\CGCtr{s_3}{s_2}{k_2}{s_1}{k_1}$
near $k_2=-s_2+i\frac{Q}{2}$ one needs to deform the
contour in \rf{CGCFT} into a small circle around the pole at $s=-R_3$ 
together with a contour separating that pole from all the other 
poles in the {\it lower} half plane. 

It is then straightforward to calculate
the values of the corresponding residues from \rf{CGCFT}. 
\end{proof}

\appendix
\section{The Clebsch-Gordan coefficients for the
modular double}

\setcounter{equation}{0}
\begin{defn}
Define a distributional kernel 
$\CGC{s_3}{x_3}{s_2}{x_2}{s_1}{x_1}$
(the ``Clebsch-Gordan coefficients'') by an expression of the form
\begin{equation}
\CGC{s_3}{x_3}{s_2}{x_2}{s_1}{x_1}\;\equiv
\; \lim_{\ep\downarrow 0}\; \CGC{s_3}{x_3}{s_2}{x_2}{s_1}{x_1}_{\ep}^{},
\end{equation}
where the meromorphic function
$\CGC{s_3}{x_3}{s_2}{x_2}{s_1}{x_1}_{\ep}^{}$ is defined as
\begin{equation}\begin{aligned}\label{Clebsch}
\CGC{-s_3}{x_3}{s_2}{x_2}{s_1}{x_1}_{\ep}
\;= \; & e^{-\frac{\pi i}{2}(h_{s_3}-h_{s_2}-h_{s_1})}\\
& \qquad\ti
D_b(\si_{32};y_{32}+i\ep)D_b(\si_{31};y_{31}+i\ep)D_b(\si_{21};y_{21}+i\ep), 
\end{aligned}\end{equation}
$h_{s}=s^2+\frac{1}{4}Q^2$, the distribution $D_b(\si;y)$ is 
defined in terms of the 
function $w_b(y)$ as 
\begin{equation}
D_b(\si;y)=\frac{w_b(y-\frac{i}{2}Q)}{w_b(y+s)},
\end{equation}
and the coefficients $y_{ji}$, $\be_{ji}$, $j>i\in\{1,2,3\}$ are given by
\begin{equation} \label{CGcoeff}
\begin{aligned}
y_{32}=& x_2-x_3+\fr{1}{2}(s_3+s_2)\\
y_{31}=& x_3-x_1+\fr{1}{2}(s_3+s_1)\\
y_{21}=& x_2-x_1+\fr{1}{2}(s_2+s_1-2s_3)
\end{aligned}
\qquad \quad \begin{aligned}
\si_{32}=& s_1-s_2-s_3\\
\si_{31}=& s_2-s_3-s_1\\
\si_{21}=& s_3-s_2-s_1.
\end{aligned}
\end{equation}
\end{defn}
It is often useful to consider the Fourier-transform
of the b-Clebsch-Gordan symbols defined by
\begin{equation}\label{FT_CGC}
\CGC{s_3}{-k_3}{s_2}{k_2}{s_1}{k_1}\;=\;
\int_{\BR}dx_3dx_2dx_1\;e^{2\pi i \sum_{l=1}^3 k_lx_l}
\CGC{s_3}{x_3}{s_2}{x_2}{s_1}{x_1}
\end{equation}
The distribution $\CGC{s_3}{k_3}{s_2}{k_2}{s_1}{k_1}$
can be factorized as
\[
\CGC{s_3}{k_3}{s_2}{k_2}{s_1}{k_1}\;=\;
\de(k_1+k_2-k_3)\CGCtr{s_3}{s_2}{k_2}{s_1}{k_1},
\]
where $\CGCtr{s_3}{s_2}{k_2}{s_1}{k_1}$
is a meromorphic function. 
A straightforward calculation using \cite[Lemma 15]{PT} yields the 
following expression:
\begin{equation}\label{CGCFT}\begin{aligned}
\CGCtr{-s_3}{s_2}{k_2}{s_1}{k_1}\;=\;
\frac{e^{-\frac{\pi i}{2}\be_{21}\be}e^{\frac{Q}{2}\pi(k_1-k_2)}e^{\pi i(k_1s_2-k_2s_1)}}{w_b(\si_{32})
w_b(\si_{31})w_b(\si_{21})}
\int\limits_{\BR+i0}ds\;e^{-\pi  s\be}
\prod_{l=1}^3\frac{w_b(s+R_l)}{w_b(s+S_l)},
\end{aligned}
\end{equation}
where we used the abbreviations $\be=\fr{Q}{2}+i(s_1+s_2+s_3)$,
$\be_{21}=\fr{Q}{2}+i(s_1+s_2-s_3)$ and
\begin{equation}  
\begin{aligned} R_1=&-s_2+k_2, \\
        R_2=& -s_1-k_1, \\
        R_3=& s_3-s_2-s_1,
\end{aligned}\qquad
\begin{aligned} 
        S_1=&  i\fr{Q}{2}+R_1-\si_{32},\\
        S_2=&  i\fr{Q}{2}+R_2-\si_{31},\\
        S_3=&  i\fr{Q}{2},
\end{aligned}
\end{equation}
The analytic properties of $\CGC{s_3}{k_3}{s_2}{k_2}{s_1}{k_1}$
can be summarized as follows:
\begin{lem}\label{CGpoles} (Lemma 20 in \cite{PT}) 
$\CGC{s_3}{k_3}{s_2}{k_2}{s_1}{k_1}$ depends meromorphically on 
all of its arguments, with poles at
$\pm ik_i=\frac{Q}{2}+is_i+nb+mb^{-1}$, $n,m\in\BZ^{\geq 0}$, $i=1,2,3$ only.
\end{lem}

\section{Quantum exponential function
 and $b$-binomial coefficient}
\setcounter{equation}{0}
The definition (\ref{gG}) and the property
(\ref{Gfunrel}) imply that the function $g_b(x)$ 
obeys the following functional equation
\begin{equation}\label{gfun}
  g_b(q x) \,=\, (1+x)^{-1} \, g_b(q^{-1}x) \,.
\end{equation}
For a pair of Weyl-type variables, $uv=q^2 vu$, 
a consequence of (\ref{gfun}) is 
\begin{equation}\label{upv}
  u+v \,=\, g_b(q u^{-1}v) \ u \ 
    \bigl(g_b(q u^{-1}v)\bigr)^{-1}   \,=\, 
 \bigl(g_b(q uv^{-1})\bigr)^{-1} \ v \ g_b(q uv^{-1}) \,.
\end{equation}
It is now obvious that
\begin{equation}\label{guv}
  g_b(u+v) = g_b(q u^{-1}v) \, g_b(u) \, 
    \bigl(g_b(q u^{-1}v)\bigr)^{-1}   
 = \bigl(g_b(q uv^{-1})\bigr)^{-1} \, 
 g_b(v) \, g_b(q uv^{-1}) \,.
\end{equation}
These relations allow us to prove the equivalence
of (\ref{qexp}) and (\ref{pent}) stated in Lemma~\ref{Qexp}.
For instance, let us show that (\ref{pent}) together
with the first relation in (\ref{guv}) implies~(\ref{qexp}).
Introduce $V\equiv q u^{-1}v$. Notice that
$V = e^{b\tilde{B}}$ where $\tilde{B}=B-A$ so that
$[A,\tilde{B}]= 2\pi i$. Then we have
\begin{equation}\nonumber
 g_b(u+v) \={guv}  g_b(V) \, g_b(u) \, 
  \bigl(g_b(V)\bigr)^{-1}   \={pent}  
 g_b(u) \, g_b(q ^{-1} u V) = g_b(u) \, g_b(v) \,. 
\end{equation}
The inverse implication, (\ref{pent}) $\Rightarrow$
(\ref{qexp}), is proven similarly.

Now we want to prove~(\ref{qexp}). First, we
represent $(u+v)^{it}$ in an integral form:
\begin{equation}\nonumber
\begin{aligned}
 (u+&v)^{it} \={upv} g_b(q u^{-1}v) \ u^{it} \ 
    \bigl(g_b(q u^{-1}v)\bigr)^{-1}  \\
 & \={GG1} b^2 \! \int\nolimits_\BR d\tau_1 d\tau_2 \, 
 e^{\pi b Q (\tau_1 -\tau_2)} \, 
 \frac{G_b(-ib\tau_1)}{G_b(Q+ib\tau_2)} \,
 (q u^{-1} v)^{i\tau_1} \, u^{it} \,
 (q u^{-1} v)^{i\tau_2}  \\
 &= b^2 \! \int\nolimits_\BR d\tau_1 d\tau_2 \, 
 e^{\pi b Q (\tau_1 -\tau_2) -i\pi b^2 (\tau_1+\tau_2)^2 
 + 2 i \pi b^2 t \tau_1} \, 
 \frac{G_b(-ib\tau_1)}{G_b(Q+ib\tau_2)} \,
 u^{i(t-\tau_1 -\tau_2)} \, v^{i(\tau_1+\tau_2)} \\
 &= b^2 \! \int\nolimits_\BR d\tau \, d\tau_2 \, 
 e^{\pi b \tau (Q +2 i b t) -i\pi b^2 \tau^2 
 - 2 \pi b \tau_2 (Q + i b t)} \, 
 \frac{G_b(ib\tau_2-ib\tau)}{G_b(Q+ib\tau_2)} \,
 u^{i(t-\tau)} \, v^{i\tau}
 \label{bpow}
\end{aligned}
\end{equation}
where we introduced $\tau\equiv \tau_1+\tau_2$. 
Computing the integral over $\tau_2$ with the help 
of~(\ref{bbeta}), we derive an analogue of the binomial formula:
\begin{equation}\label{bbin}
 (u+v)^{it} = b \int\nolimits_{\BR +i 0} d\tau \
 \BBC{t}{\tau} \, u^{i(t-\tau)} \, v^{i\tau} \,,
\end{equation}
where the $b$-binomial coefficient is given by
\begin{equation}\label{bbc}
 \BBC{t}{\tau}
 =  \frac{e^{ 2\pi i b^2 \tau (t - \tau)} \, G_b(Q+i b t)}%
 {G_b(Q+ ibt - ib \tau) \, G_b(Q + ib \tau) } \,.
\end{equation}
We see that the function $G_b$ is a $b$-analogue of the
factorial. The $b$-binomial coefficients satisfy
the $q$-Pascal identity:
\begin{equation}\label{bpas}
 \BBC{t-i}{\tau} = q^{-2 i\tau} \, \BBC{t}{\tau}
 +  \BBC{t}{\tau+i}
 =  \BBC{t}{\tau} + q^{2i(\tau -t+i)} \, \BBC{t}{\tau+i}
\end{equation}
which can be easily verified with the help of~(\ref{Gfunrel}).

Using the $b$-binomial coefficients and the integral
representation (\ref{GG1}) of $g_b(x)$, we derive the
quantum exponential relation:
\begin{equation}\nonumber
\begin{aligned}
 g_b(u+v) & \={GG1}  b\! \int\nolimits_\BR dt \, (u+v)^{it} \,
 \frac{ e^{-\pi i b^2 t^2}}{G_b(Q+ibt)} 
 \={bbin}  b^2 \!\! \int\nolimits_\BR dt \, d\tau \,
 \BBC{t}{\tau} \,
 \frac{ e^{-\pi i b^2 t^2}}{G_b(Q+ibt)} 
 u^{i(t-\tau)} v^{i\tau} \\
 & \={bbc}  b^2 \! \int\nolimits_\BR dt \, d\tau \,
  \frac{ e^{-\pi i b^2 (t-\tau)^2 -\pi i b^2 \tau^2 }}%
 {G_b(Q+ib(t-\tau)) G_b(Q+ib\tau)} \, u^{i(t-\tau)} v^{i\tau} \\
 & = b \! \int\nolimits_\BR dT \,
  \frac{ e^{-\pi i b^2 T^2 }}{G_b(Q+ibT)} \, u^{iT} \ \ 
 b \! \int\nolimits_\BR d\tau \,
 \frac{ e^{-\pi i b^2 \tau^2} }{G_b(Q+ib\tau)} \, v^{i\tau}
 \={GG1}  g_b(u) \, g_b(v) \,.
\end{aligned} 
\end{equation}
This completes the proof of Lemma~\ref{Qexp}.

\begin{rem}
After this manuscript was written we were informed that
a different proof of the quantum exponential relation  
and of relation (\ref{qF}) is given in~\cite{V}.
\end{rem}

\subsection*{Acknowledgments} 
We are grateful to L.Faddeev for useful comments.
We thank A.Volkov for providing us with a draft   
version of [V] prior to its publication.
A.B.\ was supported by
Alexander von Humboldt Foundation. 
J.T.\ was supported by DFG SFB 288.
A part of this
work was carried out during A.B.'s visit to the
Department of Mathematics, University of York.

\end{document}